\documentclass{amsart}

\newtheorem{theorem}{Theorem}[section]
\newtheorem{lemma}[theorem]{Lemma}
\newtheorem{proposition}[theorem]{Proposition}

\newtheorem{corollary}[theorem]{Corollary}

\theoremstyle{definition}

\newtheorem{example}[theorem]{Example}

\theoremstyle{remark}
\newtheorem{remark}[theorem]{Remark}

\numberwithin{equation}{section}

\begin{document}

\title[Convex-transitive characterizations of Hilbert spaces]
 {Convex-transitive characterizations of Hilbert spaces} 

\author{Jarno Talponen}
\address{Jarno Talponen, 
University of Helsinki, Department of Mathematics and Statistics, Box 68, 
(Gustaf H\"{a}llstr\"{o}minkatu 2b) FI-00014 University of Helsinki, Finland}
\email{talponen@cc.helsinki.fi}

\date{\today}


\newcommand{\no}{\addtocounter{equation}{1}(\arabic{section}.\arabic{equation})}
\renewcommand{\dim}{\mathrm{dim}}
\newcommand{\codim}{\mathrm{codim}}
\newcommand{\Mid}{\mathrm{Mid}}
\renewcommand{\span}{\mathrm{span}}
\newcommand{\diam}{\mathrm{diam}} 
\newcommand{\dist}{\mathrm{dist}}
\newcommand{\conv}{\mathrm{conv}}
\newcommand{\card}{\mathrm{card}}
\newcommand{\E}{\mathbf{E}}
\renewcommand{\P}{\mathbf{P}}
\newcommand{\R}{\mathbb{R}}
\newcommand{\N}{\mathbb{N}}
\newcommand{\C}{\mathbb{C}}
\newcommand{\K}{\mathbb{K}}
\newcommand{\Q}{\mathbb{Q}}
\newcommand{\X}{\boldsymbol{X}}
\renewcommand{\H}{\boldsymbol{H}}
\newcommand{\Y}{\boldsymbol{Y}}
\newcommand{\Z}{\boldsymbol{Z}}
\newcommand{\B}{\mathbf{B}}
\newcommand{\I}{\mathbf{I}}
\newcommand{\sign}{\mathrm{sign}}
\renewcommand{\S}{\mathbf{S}}
\renewcommand{\b}{\boldsymbol{b}}
\newcommand{\vv}{\boldsymbol{v}}
\renewcommand{\t}{\boldsymbol{t}}
\renewcommand{\ker}{\mathrm{Ker}}
\newcommand{\T}{\mathbb{T}}
\newcommand{\el}{\ell}
\newcommand{\wexp}{\omega\mathrm{-exp}}
\newcommand{\nn}{\mathfrak{n}_{k}}
\newcommand{\mm}{\mathfrak{m}_{l}}
\maketitle

\begin{abstract}
In this paper we investigate real convex-transitive Banach spaces $\X$, which admit a $1$-dimensional bicontractive projection 
$P$ on $\X$. Various mild conditions regarding the weak topology and the geometry of the norm are provided, 
which guarantee that such an $\X$ is in fact isometrically a Hilbert space. The results obtained here can be regarded as partial 
answers to the well-known Banach-Mazur rotation problem, as well as to a question posed by B. Randrianantoanina in 2002 about 
convex-transitive spaces.  
\end{abstract}

\section{Introduction}
This work draws its motivation from the \textit{Banach-Mazur rotation problem}, which was already formulated 
in Banach's book \cite[p.242]{Ba}. We denote by $\X$ real Banach spaces, by $\B_{\X}$ the closed unit ball of $\X$
and by $\S_{\X}$ the unit sphere of $\X$. The orbit of $x$ under the group of rotations
$\mathcal{G}_{\X}=\{T|\ T\in L(\X)\ \mathrm{isometry\ onto}\}$ is given by $\mathcal{G}_{\X}(x)=\{T(x)|\ T\in \mathcal{G}_{\X}\}$.
Recall that $\X$ is \emph{transitive} if $\mathcal{G}_{\X}(x)=\S_{\X}$ for all $x\in\S_{\X}$, and 
that $\X$ is \emph{almost transitive} if $\overline{\mathcal{G}_{\X}(x)}=\S_{\X}$ for all $x\in\S_{\X}$.
An element $x\in\S_{\X}$ is called a \emph{big point} if $\overline{\conv}(\mathcal{G}_{\X}(x))=\B_{\X}$ and 
$\X$ is \emph{convex-transitive} if each $x\in\S_{\X}$ is a big point. Clearly Hilbert spaces are transitive.
Intuitively speaking, the unit ball of a transitive space can be thought of as roundish. 
The following \textit{rotation problem} remains unanswered:
\begin{flushleft}
($\mathbf{Q}_{\mathrm{BM}}$)\ \textit{Is every separable transitive Banach space $\X$ in fact isometrically $\ell^{2}$?}
\end{flushleft}
Recall that there exists non-separable transitive spaces $L^{p}(\Omega,\mu)$ for $1\leq p<\infty$, see \cite{Rol}.
Extensive surveys of the wide body of results generated around the problem ($\mathbf{Q}_{\mathrm{BM}}$) are found in \cite{BR2} 
and \cite{Ca1}. 
 
Some partial answers to the rotation problem can be given by imposing a variety of additional conditions.
In this connection it is natural to require some weak substitutes for orthogonality. 
For instance, each of the following isometric conditions \addtocounter{equation}{1} $(\arabic{section}.\arabic{equation})$-\addtocounter{equation}{2}$(\arabic{section}.\arabic{equation})$\addtocounter{equation}{-3} characterizes Hilbert spaces isometrically:
\begin{flushleft}
\no\  $\X$ is convex-transitive and for some $1\leq p\leq\infty,$ there is a $1$-dimensional subspace $L\subset\X$ such that 
$\X=L\oplus\Y$, and $||(z,y)||=||(||z||_{L},||y||_{\Y})||_{\ell^{p}_{2}}$ for all $z\in L$ and $y\in\Y$ (see \cite{BR1}).\\
\newcounter{sect_a1}
\setcounter{sect_a1}{\value{section}}
\newcounter{eq_a1}
\setcounter{eq_a1}{\value{equation}}
\no\ $\X$ is almost transitive and there is a $1$-dimensional subspace $L\subset\X$ such that 
$\X=L\oplus\Y$, where $||(z,y)||_{\X}=||(-z,y)||_{\X}$ for all $z\in L,\ y\in\Y$ (see \cite{SZ}).\\
\newcounter{sect_a2}
\setcounter{sect_a2}{\value{section}}
\newcounter{eq_a2}
\setcounter{eq_a2}{\value{equation}}
\no\ $\X$ is almost transitive and there is a $1$-dimensional subspace $L\subset\X$ such that $\X=L\oplus\Y$, 
where $\Y$ is $1$-complemented (see \cite{Ra2}).
\newcounter{sect_a3}
\setcounter{sect_a3}{\value{section}}
\newcounter{eq_a3}
\setcounter{eq_a3}{\value{equation}}
\end{flushleft}

Note that the preceding conditions about decompositions of orthogonal type appear in weakening order. 
In characterization $(\arabic{sect_a3}.\arabic{eq_a3})$ there is, as it turns out (see section 1.1), a linear closest point selection 
$P\colon \X\rightarrow L$. 
The question was raised by B. Randrianantoanina in \cite{Ra2} whether the characterization 
$(\arabic{sect_a3}.\arabic{eq_a3})$ can be extended to the convex-transitive setting:
\begin{flushleft}
($\mathbf{Q}_{\mathrm{R}}$)\ \textit{Are all convex-transitive spaces $\X$ which contain a $1$-codimensional, $1$-complemented 
subspace $\Y\subset\X$ in fact isometrically Hilbert spaces?}
\end{flushleft}
In fact, it is not even known if each \textit{real} convex-transitive Banach space with an \textit{isometric reflection vector}, 
as in $(\arabic{sect_a2}.\arabic{eq_a2})$, is isometrically a Hilbert space. Since the decompositions above appear in weakening order, 
any generalization of $(\arabic{sect_a3}.\arabic{eq_a3})$ also applies to $(\arabic{sect_a2}.\arabic{eq_a2})$ as well. 
The main part of this paper deals with question $(\mathbf{Q}_{\mathrm{R}})$.

It turns out that the answer to question ($\mathbf{Q}_{\mathrm{R}}$) is \emph{affirmative} in a very wide class of spaces. 
Compared to the original paper \cite{Ra2} a novel approach is required to tackle ($\mathbf{Q}_{\mathrm{R}}$).
In fact, our approach requires some mild assumptions about the interplay of the weak topology and the geometry of the norm
(in some cases these additional assumptions can not be dispensed with).
This interplay is termed here the \emph{weak norm geometry}. It involves, roughly speaking, the properties of 
$(\S_{\X},\omega)$ and is comprised of e.g. weak local uniform rotundity and G\^{a}teaux smoothness. For related geometric results based 
on assumptions in terms of the \emph{norm} topology, see \cite{BR1,BR3}.

On the other hand, a suitable control of the weak norm geometry appears to be the natural approach in conjunction with the 
convex-transitivity. This is so because control of the norm geometry (in the sense of \emph{norm} topology) in fact often reduces 
the convex-transitive case to the almost transitive one, see e.g. \cite[p.51-54]{BR2}. 
 
The following result (Theorem \ref{th:wlurpoint}) related to question $(\mathbf{Q}_{\mathrm{R}})$ is essential here:\\
\textit{Suppose that $\X$ is a convex-transitive Banach space, which admits a $\omega\mathrm{-LUR}$ point $x\in\S_{\X}$. Assume further
that $\X$ admits a $1$-dimensional bicontractive projection $P\colon \X\rightarrow L$. Then $\X$ is isometrically a Hilbert space.}

Theorem \ref{bpt}, which is another main result, illustrates the fact that if the assumption about the weak norm geometry 
is strengthened, then the requirement of convex-transitivity of $\X$ can be relaxed to one of $\X$ admitting a big point.

This paper is organized as follows. Sections 3 and 5 contain the main results, which characterize Hilbert spaces.
In section 3 the characterizations are formulated in terms of big points and in section 5 by assuming convex-transitivity of the space.
In Section 2 some common technical components of the subsequent proofs are extracted and formulated in rather general form. 
We apply, roughly speaking, the following three main strategies in attacking the problem $(\mathbf{Q}_{\mathrm{R}})$. 
In order to check that a given projection $P\colon \X \rightarrow [u]$ has the desired properties (such as bicontractiveness) we 
reconstruct this projection by suitable approximation with operators having such properties. Lemma \ref{lm2} makes use of the 
$\omega^{\ast}$-compactness of $\B_{\X^{\ast}}$ and the density of the smooth points $y\in \S_{\X}$. 
This lemma will be a powerful tool for approximating $1$-dimensional operators by large-dimensional operators. 

In order to apply Lemma \ref{lm2}, which requires smoothness, we will reduce our setting (in Theorem \ref{extpro}) to a suitable
separable subspace of the given space $\X$. This result is proved by applying a back-and-forth type recursion of countable length. 

One common feature of the main results in sections 3 and 5, apart from the above mentioned approximation, is that the weak norm 
geometry acts as a kind of glue, which conveys the orbits in close proximity in the sense of the weak topology. 

In section 4 we introduce a concept of atoms for general Banach spaces, which is motivated by atoms in Banach lattices.
Both these concepts of atoms are included in the more general class of isometric reflection vectors.
We show that a Banach space admitting an atom, which is simultaneously a big point, is in fact isometrically 
$\ell^{1}(\Gamma)$ for a suitable set $\Gamma$. 

Finally, we point out that the characterization $(\arabic{sect_a3}.\arabic{eq_a3})$ has been recently extended to an almost isometric 
setting in \cite{asypri}. 
\subsection{Preliminaries}
In what follows $\X$ is a real Banach space and $\X^{\ast}$ is its dual space.
The weak topology of $\X$ and the weak-star topology of $\X^{\ast}$ are denoted by $\omega$ and by $\omega^{\ast}$, respectively. 
Here $\N=\{1,2,3,\ldots\}$ but $0\in\omega_{0}$, the first infinite ordinal.

We denote by $L(\X,\Y)$ the space of continuous linear operators $T\colon \X\rightarrow\Y$ and
$L(\X)=L(\X,\X)$ for short. An \emph{isomorphism} $T\in L(\X)$ is called an \emph{automorphism}, denoted $T\in\mathrm{Aut}(\X)$. 
If $T\in\mathrm{Aut}(\X)$ is an isometry then we say that $T$ is a \emph{rotation}.
We denote by $\mathcal{G}=\mathcal{G}_{\X}=\{T\colon \X\rightarrow \X\ |\ T\ \mathrm{is\ a\ rotation}\}$ 
the \emph{group of rotations} of $\X$. The composition of the maps acts as the group operation and the identity map 
$\I\colon \X\rightarrow\X$ is the neutral element. 

Recall that a point $x\in \S_{\X}$ is \textit{weakly exposed} (in $\B_{\X}$) by $f\in\S_{\X^{\ast}}$  if for each sequence
$(x_{n})\subset\B_{\X}$ such that $f(x_{n})\rightarrow 1$ as $n\rightarrow\infty$ it holds that 
$x_{n}\stackrel{\omega}{\longrightarrow}x$ as $n\rightarrow\infty$. Denote the weakly exposed points in $\B_{\X}$ by $\wexp(\B_{\X})$.
If $\tau$ is a given locally convex topology on $\X$ then $x\in\S_{\X}$ is said to be $\tau$-LUR (short for locally uniformly rotund) 
if for any sequence $(x_{n})\subset\B_{\X}$ such that $||x+x_{n}||\rightarrow 2$ as $n\rightarrow\infty$ it holds that 
$x_{n}\stackrel{\tau}{\longrightarrow}x$ as $n\rightarrow\infty$.
A space $\X$ is called $\tau$-LUR if each point $x\in\S_{\X}$ is $\tau$-LUR.  
An \emph{open slice} has the form $S_{f,\alpha}=\{x\in\B_{\X}|f(x)>\alpha\}$, where $f\in\X^{\ast}$ and $\alpha\in\R$.
We will frequently apply the following classical results (see \cite[p.96,92]{HHZ}). 
\begin{lemma}[\v Smulyan]\label{smulyan1}
The following conditions \addtocounter{equation}{1} $(\arabic{section}.\arabic{equation})$-\addtocounter{equation}{2}$(\arabic{section}.\arabic{equation})$\addtocounter{equation}{-3} are equivalent for any Banach space $\X$:
\begin{enumerate}
\item[\no]{The norm $||\cdot||$ of $\X$ is G\^{a}teaux smooth at $x\in \S_{\X}$}
\newcounter{sect_a4}
\setcounter{sect_a4}{\value{section}}
\newcounter{eq_a4}
\setcounter{eq_a4}{\value{equation}}
\item[\no]{For all $(f_{n}),(g_{n})\subset \S_{\X^{\ast}}$ such that 
$\lim_{n\rightarrow\infty}f_{n}(x)=\lim_{n\rightarrow\infty}g_{n}(x)=1$ it holds that $f_{n}-g_{n}\stackrel{\omega^{\ast}}{\longrightarrow} 0\ \mathrm{as}\ n\rightarrow\infty$.}
\item[\no]{There is unique $f\in \S_{\X^{\ast}}$ such that $f(x)=1$.}
\newcounter{sect_a5}
\setcounter{sect_a5}{\value{section}}
\newcounter{eq_a5}
\setcounter{eq_a5}{\value{equation}}
\end{enumerate}
\end{lemma}
\begin{theorem}[(Mazur)]\label{classmazur}
Let $\X$ be a separable Banach space. Then the set of G\^{a}teaux smooth points $x\in \S_{\X}$ is a dense $G_{\delta}$-set of $\S_{X}$.
\end{theorem}

Observe that if $x\in\S_{\X}$ is smooth and $f\in \S_{\X^{\ast}}$ satisfies $f(x)=1$ then the family 
\begin{equation}\label{eq: neighbourhood basis}
\{g\in \B_{\X^{\ast}}|g(x)>1-\frac{1}{n}\}\subset \B_{\X^{\ast}},\quad n\in\N
\end{equation}
defines a $\omega^{\ast}$-neighbourhood basis for $f$ relative to $(B_{\X^{\ast}},\omega^{\ast})$. Indeed, let $U\subset \B_{\X^{\ast}}$ be a $\omega^{\ast}$-open neighbourhood for $f$. Then there is $k\in \N$ such that $\{g\in \B_{\X^{\ast}}|g(x)>1-\frac{1}{n}\}\subset U$ 
for all $n\geq k$. In fact, otherwise one could pick a subsequence $(n_{l})_{l}\subset \N$ together with 
$g_{l}\in \{g\in \B_{\X^{\ast}}|g(x)>1-\frac{1}{n_{l}}\}\setminus U$ for $l\in\N$. On the other hand, according to the \v Smulyan 
lemma $g_{l}\stackrel{\omega^{\ast}}{\longrightarrow}f$ as $l\rightarrow\infty$. Since $U$ is an open neighbourhood of $f$ this 
contradicts the fact that $\{g_{l}|l\in\N\}\cap U=\emptyset$.

Next we will mention some facts about the closest point projections and Flinn pairs, which are the essential concepts in this article 
together with the transitivity conditions. In what follows we let $f\otimes x\colon \X\rightarrow [x]$ be the mapping
$y\longmapsto f(y)x$, where $f\in\X^{\ast},\ x\in\X$. 
Following \cite{KR} and \cite{Ra2} we call $(x,f)\in \S_{\X}\times\X^{\ast}$ a \textit{Flinn pair} if 
$||\I-f\otimes x||=1$. In such a case the $x$ above is called a \emph{Flinn element}. Hence Randrianantoanina's 
question $(\mathbf{Q}_{\mathrm{R}})$ from the introduction involves the existence of Flinn elements
in the convex-transitive setting. See \cite{Ra2} for connections between Flinn pairs and numerically positive 
operators and \cite{Ra1} for related results.

As an example, consider a linear projection $P\colon \X\rightarrow [x]$, which is an 
\emph{isometric reflection projection}, that is, $2P-\I\in \mathcal{G}_{\X}$. 
Then it is easy to see that $P=f\otimes x$ for some Flinn pair $(x,f)\in \S_{\X}\times \S_{X^{\ast}}$.

Recall that for a subset $A\subset \X$ a map $P_{A}\colon \X\rightarrow A$
is called a closest point selection if $\dist(x,A)=||x-P_{A}(x)||$ for all $x\in\X$.
In \cite[Lemma 13.1]{Am} it is pointed out that for a closed linear subspace $A\subset \X$ a linear projection 
$P_{A}\colon \X \rightarrow A$ is a closest point selection if and only if $||\I-P_{A}||=1$. 
Consequently $(x,f)$ is a Flinn pair if and only if $f\otimes x$ is a linear closest point selection to $[x]$. 
In other words, Flinn elements $x$ are exactly those points of $\S_{\X}$ for which there exists a linear closest point selection 
$P_{[x]}\colon \X \rightarrow [x]$. This allows us to use the terminology '$(x,f)$ is a Flinn pair' and
'$f\otimes x$ is a closest point selection' interchangeably. 

The following characterization of Hilbert spaces will be applied frequently, (see \cite[p.108]{Am}).
\begin{proposition}\label{char}
A Banach space $\X,\ \dim(\X)\geq 3,$ is isometrically a Hilbert space if and only if for each $1$-dimensional subspace 
$L\subset \X$ there exists a linear projection $P\colon \X\rightarrow L$ such that $||\I-P||=1$.
\end{proposition}

In the attempts to answer the question $(\mathbf{Q}_{\mathrm{R}})$ we will often assume additionally that above $||P||=1$.
One can argue that this is not too restrictive an assumption, since any orthogonal projection 
$P\colon \H\rightarrow [x]$, or more generally, any isometric reflection projection $P\colon\X\rightarrow [x]$ 
satisfies $||P||=||\I-P||=1$. In this event $P$ is called a \emph{bicontractive} projection.

\section{Machinery}
The following lemma might be known in some form but a proof is included here for convenience. 
\begin{lemma}\label{lm}
Let $\X$ be a Banach space. Assume that $g\in \S_{\X^{\ast}}$ is a support functional of $x\in \S_{\X}$. If the sequence 
$(x_{n})\subset \B_{\X}$ satisfies $x\in\overline{\conv}(\{x_{n}|n\in \N\})$, then: 
\begin{enumerate}
\item[(i)]{For each $\epsilon>0$ it holds that $x\in\overline{\conv}(\{x_{n}|g(x_{n})>1-\epsilon\})$.}
\item[(ii)]{There are finite sets $K_{l}\subset \N,\ l\in \N,$ such that for each $l\in \N$\\
$\inf\{g(x_{n})|n\in K_{l}\}\geq 1-\frac{1}{l}\ \mathrm{and}\ \dist(x,\conv(\{x_{n}|n\in K_{l}\}))<\frac{1}{l}$.}
\item[(iii)]{There is a sequence $(y_{j})\subset \{x_{n}|n\in \N\}$ such that $g(y_{j})\rightarrow 1$ as $j\rightarrow \infty$
and \mbox{$x\in \overline{\conv}(\{y_{j}|j\geq l\})$} for each $l\in \N$.}
\end{enumerate}
\end{lemma}   
\begin{proof}
We first prove the main claim (i). Fix $\epsilon>0$.
Let $(c^{(k)}_{n})_{n}\in \S_{\ell_{+}^{1}}\cap c_{00}$ be a sequence of non-negative convex coefficients for each $k\in\N$ such that 
$\sum_{n}c_{n}^{(k)}x_{n}\rightarrow x$ as $k\rightarrow\infty$. 
Write $P_{\epsilon}^{(k)}=\sum_{n:g(x_{n})\geq 1-\epsilon}c_{n}^{(k)}$ and $Q_{\epsilon}^{(k)}=1-P_{\epsilon}^{(k)}$ for $k\in\N$. Then
\begin{equation}\label{eq: Cnk}
\begin{array}{rcl}
\sum_{n}c_{n}^{(k)}x_{n}&=&P_{\epsilon}^{(k)}\frac{\sum_{n:g(x_{n})\geq 1-\epsilon}c_{n}^{(k)}x_{n}}{P_{\epsilon}^{(k)}}+Q_{\epsilon}^{(k)}\frac{\sum_{n:g(x_{n})<1-\epsilon}c_{n}^{(k)}x_{n}}{Q_{\epsilon}^{(k)}},
\end{array}
\end{equation}
where we apply the convention $\frac{0}{0}=0$. Observe that above 
\[\sum_{n:g(x_{n})<1-\epsilon}\frac{c_{n}^{(k)}}{Q_{\epsilon}^{(k)}}x_{n}\ \mathrm{and}\ \sum_{n:g(x_{n})\geq 1-\epsilon}\frac{c_{n}^{(k)} }{P_{\epsilon}^{(k)}}x_{n}\]
are members of $\conv(\{x_{n}|n\in\N\})$ for all $k\in\N$. Note that 
\[Q_{\epsilon}^{(k)}(1-\epsilon)+P_{\epsilon}^{(k)}\geq \sum_{n}c_{n}^{(k)}g(x_{n})
=g(\sum_{n}c_{n}^{(k)}x_{n})\stackrel{k\rightarrow\infty}{\longrightarrow}g(x)=1.\]
Since $Q_{\epsilon}^{(k)}+P_{\epsilon}^{(k)}=1$ for all $k\in\N$, 
we obtain that $Q_{\epsilon}^{(k)}\rightarrow 0$ as $k\rightarrow \infty$.

Hence we obtain that the norm of the second term in \eqref{eq: Cnk} tends to $0$ as $k$ goes to infinity. Consequently
\[\sum_{n}c_{n}^{(k)}x_{n}-\sum_{n:g(x_{n})\geq 1-\epsilon}\frac{c_{n}^{(k)}}{P_{\epsilon}^{(k)}}x_{n}\longrightarrow 0\] 
as $k\rightarrow\infty$. This means that $x\in\overline{\conv}(\{x_{n}|g(x_{n})\geq 1-\epsilon\})$, which is claim (i). 
By applying this case it is then easy to select suitable indices for the claims (ii)-(iii).
\end{proof}
We would like to highlight the following case, which will subsequently become relevant.
Suppose that $\{(c^{(k)}_{n})_{n}|k\in\N\}\subset \S_{\ell^{1}_{+}}\cap c_{00}$, $(V,\tau)$ is a topological vector space, $x\in V$ and 
$(x_{n})\subset V$ is a sequence. Let us consider the case that 
\begin{equation}\label{eq: limU}
\lim_{k\rightarrow\infty}\sum_{n: x_{n}\in U}c_{n}^{(k)}=1\quad \mathrm{whenever}\ U\subset V\ \mathrm{is\ a}\ \tau\mathrm{-open\ neighbourhood\ of}\ x.
\end{equation}
The following condition is a reformulation of \eqref{eq: limU}: Whenever $(U_{j})\subset V$ is a sequence
of $\tau$-open neighbourhoods of $x$, then there exists a non-decreasing map $\alpha\colon \N\rightarrow \N$ such that 
\begin{equation}\label{eq: equivalent}
\lim_{i\rightarrow\infty}\alpha(i)=\infty\ \mathrm{and}\ \lim_{k\rightarrow\infty}\sum_{n:x_{n}\in \bigcap_{j\leq\alpha(k)}U_{j}}c_{n}^{(k)}=1.
\end{equation}
In fact, put 
\[\alpha(i)=\max\left\{\alpha\in\N|\alpha\leq i\ \mathrm{and}\ \sum_{n:x_{n}\notin\bigcap_{j\leq\alpha}U_{j}}c_{n}^{(k)}\leq 2^{-\alpha+1}\ \mathrm{for\ all}\ k\geq i\right\}.\]
Observe that $\alpha(i)$ is defined for all $i\in\N$ since $2^{-\alpha +1}=1$ for $\alpha=1$.
Since the condition on the right hand side weakens as $i$ increases, we obtain that $\alpha(i)$ is non-decreasing.
If $\lim_{i\rightarrow\infty}\alpha(i)=m<\infty$, then 
\[\limsup_{k\rightarrow\infty}\sum_{n:x_{n}\notin\bigcap_{j\leq m+1}U_{j}}c_{n}^{(k)}\geq 2^{-m},\] 
which contradicts \eqref{eq: limU}, since $\bigcap_{j\leq m+1}U_{j}$ is a $\tau$-open neighbourhood of $x$.

\begin{lemma}\label{fact}
Let $\X$ be a Banach space, $(T_{j})_{j}\subset\mathcal{G}_{\X}$ and $x,y,z\in \X$. Then 
\[\dist(y, \conv(\{T_{j}(z)|j=1,\ldots,n\}))\leq \dist(y,\conv(\{T_{j}(x)|j=1,\ldots,n\}))+||x-z||.\]
Alternatively, if $(c_{n})_{n}\in \S_{\ell^{1}_{+}}\cap c_{00}$, then
\[||y-\sum_{n}c_{n}T_{n}(z)||\leq ||y-\sum_{n}c_{n}T_{n}(x)||+||x-z||.\]
\begin{proof}
We obtain 
\[||y-\sum_{n}c_{n}T_{n}(z)||\leq ||y-\sum_{n} c_{n}T_{n}(x)||+||\sum_{n} c_{n}(T_{n}(x)-T_{n}(z))||,\]
where $||\sum_{n}c_{n}(T_{n}(x)-T_{n}(z))||\leq\sum_{n}c_{n}||T_{n}(x-z)||=||x-z||$.
\end{proof}
\end{lemma}
 
The following lemma is the main technical tool applied below in the main theorems, for example in Theorem \ref{bpt}. 
It seems likely that it could have further applications combined with Theorem \ref{extpro} and
the techniques used in the proofs of Theorems \ref{bpt}, \ref{ct} and \ref{th:wlurpoint}.
\begin{lemma}\label{lm2}
Let $\X$ be a Banach space and let $(x_{n},f_{n})\in \S_{\X}\times \S_{\X^{\ast}},\ n\in \N,\ $ be pairs
such that $f_{n}(x_{n})=1$ for each $n$. Assume that $(y,g)\in \S_{\X}\times \S_{\X^{\ast}}$ satisfies $g(y)=1$ and $y$ is a smooth
point. Let $\{(c_{n}^{(k)})_{n}|k\in\N\}\subset \S_{\ell^{1}_{+}}\cap c_{00}$ be such that 
\begin{enumerate}  
\item[(i)]{$\sum_{n}c_{n}^{(k)}x_{n}\stackrel{\omega}{\longrightarrow} y$ as $k\rightarrow \infty$,}
\item[(ii)]{$\sum_{n}c_{n}^{(k)}f_{n}(y)\rightarrow 1$ as $k\rightarrow \infty$.}
\end{enumerate}
Then 
\[P_{k}\stackrel{\cdot}{=}\sum_{n}c_{n}^{(k)}f_{n}\otimes x_{n}\stackrel{\mathrm{WOT}}{\longrightarrow}g\otimes y\ \mathrm{as}\ k\rightarrow \infty,\] 
and there is $\{(d_{n}^{(l)})_{n}|l\in\N\}\subset \S_{\ell^{1}_{+}}\cap c_{00}$ such that
\[S_{l}\stackrel{\cdot}{=}\sum_{n}d_{n}^{(l)}f_{n}\otimes x_{n}\stackrel{\mathrm{SOT}}{\longrightarrow}g\otimes y\ \mathrm{as}\ l\rightarrow \infty.\] 
Moreover, if in addition $(x_{n},f_{n})$ are Flinn pairs for each $n\in \N$, then $(y,g)$ is also one.

\begin{proof}
Clearly the operators $P_{k}\colon \X\rightarrow \X$ are linear and 
\[||P_{k}||=||\sum_{n}c_{n}^{(k)}f_{n}\otimes x_{n}||\leq \sum_{n}c_{n}^{(k)}||f_{n}\otimes x_{n}||=1.\]
By assumption (ii) we have $\lim_{k\rightarrow\infty}\sum_{n}c_{n}^{(k)}(1-f_{n}(y))=0$ and 
we obtain similarly as in \eqref{eq: Cnk} that
\begin{equation}
\sum_{n:1-f_{n}(y)<1-\epsilon}c_{n}^{(k)}\longrightarrow 0\ \mathrm{as}\ k\rightarrow\infty\ \mathrm{for\ all}\ \epsilon>0.
\end{equation}
 
Thus according to \eqref{eq: equivalent} there exists a sequence $(l_{k})_{k\in\N}\subset\N$ such that $l_{k}\rightarrow\infty$ as 
$k\rightarrow\infty$ and
\begin{equation}\label{eq: Pfy}
\sum_{n:f_{n}(y)\geq 1-\frac{1}{l_{k}}}c_{n}^{(k)}\rightarrow 1\ \quad \mathrm{as}\ k\rightarrow\infty . 
\end{equation}
According to the remark related to \eqref{eq: neighbourhood basis} we obtain that
\begin{equation}\label{eq: nb}
\left(\left\{h\in \B_{\X^{\ast}}|h(y)\geq 1-\frac{1}{l_{k}}\right\}\right)_{k\in\N}\subset (\B_{\X^{\ast}},\omega^{\ast})
\end{equation}
is a neighbourhood basis for $g$, that is, for each $\omega^{\ast}$-open neighbourhood 
$U\subset \B_{\X^{\ast}}$ of $g$ there is $k$ such that $\{h\in \B_{\X^{\ast}}|h(y)\geq 1-\frac{1}{l_{k}}\}\subset U$.
By taking into account \eqref{eq: Pfy} we obtain that $\lim_{k\rightarrow\infty}\sum_{n:f_{n}\in U}c_{n}^{(k)}=1$
for any $\omega^{\ast}$-open neighbourhood $U\subset \B_{\X^{\ast}}$ of $g$. In particular, since
$f_{n}\otimes x_{n}-g\otimes x_{n}=(f_{n}-g)\otimes x_{n}$ we have 
\[\lim_{k\rightarrow\infty}\sum_{n:|(f_{n}-g)(v)|\geq\epsilon}c_{n}^{(k)}=0\ \mathrm{for\ all}\ v\in\X\ \mathrm{and}\ \epsilon>0.\]
Thus we obtain  
\[||\sum_{n}c_{n}^{(k)}(f_{n}\otimes x_{n}(v)-g\otimes x_{n}(v))||\leq \sum_{n}c_{n}^{(k)}||(f_{n}\otimes x_{n})(v)-(g\otimes x_{n})(v)||
\rightarrow 0\]
as $k\rightarrow\infty$ for all $v\in\X$. This means that 
\begin{equation}\label{eq: Efx}
\mathrm{SOT-}\lim_{k\rightarrow\infty}\sum_{n}c_{n}^{(k)}(f_{n}\otimes x_{n}-g\otimes x_{n})=0
\end{equation}
and by assumption (i) we obtain
\begin{equation}
\mathrm{WOT-}\lim_{k\rightarrow\infty}\sum_{n}c_{n}^{(k)}(g\otimes x_{n})=g\otimes(\omega\mathrm{-}\lim_{k\rightarrow\infty}\sum_{n}c_{n}^{(k)}x_{n})=g\otimes y.
\end{equation}
We conclude that $P_{k}=\sum_{n}c_{n}^{(k)}f_{n}\otimes x_{n}\stackrel{\mathrm{WOT}}{\longrightarrow}g\otimes y$ as $k\rightarrow\infty$,
which completes the first part of the lemma.
 
Since $P_{k}(y)\stackrel{\omega}{\longrightarrow}y$ in $\X$ as $k\rightarrow \infty$ we get by an application of Mazur's classical theorem 
that there is a sequence $(S_{l})\subset L(\X)$ (depending on the fixed $y$) such that $S_{l}\in \conv(\{P_{k}|k\geq l\})$ 
for each $l\in \N$ and $S_{l}(y)\stackrel{||\cdot||}{\longrightarrow} y$ as $l\rightarrow\infty$.
Write $S_{l}$ as a finite convex combination 
$S_{l}=\sum_{k\geq l} e_{k}^{(l)} P_{k},$ where $\{(e_{k}^{(l)})_{k}|l\in\N\}\subset \S_{\ell_{+}^{1}}\cap c_{00}$ 
are the corresponding coefficients. Now for each $z\in \ker(g)$ it holds by \eqref{eq: Efx} that 
$P_{k}(z)\stackrel{||\cdot||}{\longrightarrow} 0$ as $k\rightarrow \infty$. 
It follows that $\mathrm{SOT-}\lim_{l\rightarrow\infty}S_{l}=g\otimes y$. Indeed, write $v\in \X$ as $v=ay+z$, where 
$z\in \ker(g)$ and $a\in \R$. Then $S_{l}(v)=S_{l}(ay)+S_{l}(z)$, where $S_{l}(ay)\rightarrow ay$ as $l\rightarrow\infty$ and 
\[||S_{l}(z)||=||\sum_{k\geq l}e_{k}^{(l)}P_{k}(z)||\leq \sum_{k\geq l}e_{k}^{(l)}||P_{k}(z)||\rightarrow 0\ \mathrm{as}\ l\rightarrow \infty.\]
Hence $\mathrm{SOT-}\lim_{l\rightarrow\infty}S_{l}=g\otimes y$ and the corresponding convex weights $d_{n}^{(l)}\geq 0$ are given by 
$d_{n}^{(l)}=\sum_{k\geq l}e_{k}^{(l)}c_{n}^{(k)}$ for all $n,l\in\N$.
 
Finally, if $(x_{n},f_{n})$ are Flinn pairs for each $n\in\N$, then we get
\begin{eqnarray*}
||\I-g\otimes y||&\leq&\liminf_{l\rightarrow \infty} ||\I-S_{l}||=\liminf_{l\rightarrow\infty}||\sum_{n}d_{n}^{(l)}\I-\sum_{n}d_{n}^{(l)}f_{n}\otimes x_{n}||\\
		 &\leq&\liminf_{l\rightarrow\infty}\sum_{n}d_{n}^{(l)}||\I-f_{n}\otimes x_{n}||=1.
\end{eqnarray*}
This means that $||\I-g\otimes y||=1$ since $\I-g\otimes y$ is a projection and thus $(y,g)$ is a Flinn pair. 
\end{proof}
\end{lemma}

Recall that the \emph{density character} of a topological space $(T,\tau)$ is 
$$\mathrm{dens}(T)=\min\{\kappa\ \mathrm{cardinal}|\ \mathrm{there\ is\ a\ dense\ set\ } A\subset T\ \mathrm{such\ that\ } |A|=\kappa\}.$$
Recall that the first infinite cardinal $|\N|$ is denoted by $\aleph_{0}=\omega_{0}$ and
that for any infinite cardinal $\kappa$ it holds that 
$|\kappa\times\kappa|=\kappa$ (see e.g. \cite[p.162]{En}). Consequently $|\kappa^{n}|=\kappa$ for all $n<\omega_{0}$ and
$|\lambda\times \kappa|=\kappa$, where $\lambda\neq\emptyset$ is an ordinal and $\kappa\geq \lambda$ is an infinite cardinal.

Separable Banach spaces have very convenient smoothness properties, as demonstrated by Theorem \ref{classmazur}. 
In proving the main results our strategy is to reduce the non-separable case to a separable one, in order to get access to the 
applications of smoothness. The following result is a central tool in this reduction. 
\begin{theorem}\label{extpro}
Let $\X$ be a Banach space which is convex-transitive with respect to a given subgroup of isometries 
$\mathcal{G}^{0}\subset\mathcal{G}_{\X}$,
and let $C\subset \S_{\X}$ be a given norm-dense subset. Assume that $A\subset \X$ is a closed subspace with 
density character $\mathrm{dens}(A)=\kappa$. Then there exists a closed subspace $\Y\subset \X$ such that
\begin{enumerate}
\item{$A\subset \Y$}
\item{$\mathrm{dens}(\Y)=\kappa$}
\item{$\Y$ is convex-transitive with respect to the subgroup\\ 
$\mathcal{G}_{\Y}^{0}=\{T_{|\Y}|\ T\in \mathcal{G}^{0},\ T(\Y)=\Y\}\subset \mathcal{G}_{\Y}$}
\item{$\overline{C\cap \S_{\Y}}=\S_{\Y}$.}
\end{enumerate}
\end{theorem}
\begin{proof}
Clearly $\mathcal{G}_{\Y}^{0}$ above is a subgroup which at least contains the identity element $\I_{\Y}$. 
We apply a back-and-forth recursion of length $\omega_{0}$. In this proof we adopt set-theoretic notations 
as in \cite{En}.
 
Let $C\subset \S_{\X}$ be the fixed dense set. Since the density character of $A$ is $\kappa$ and $C$ is dense in $\S_{\X}$ 
we can find a set $C_{0}=\{x_{\alpha}\}_{0\leq \alpha\leq \kappa}\subset C$ such that $\S_{A}\subset\overline{C_{0}}$. 
Indeed, let $\{y_{\alpha}\}_{0\leq \alpha\leq \kappa}\subset \S_{A}$ be a dense set.
For each $\alpha$ there is a a sequence $(y_{\alpha}^{(n)})_{n<\omega_{0}}\subset C$ such that $y_{\alpha}^{(n)}\rightarrow y_{\alpha}$ as 
$n\rightarrow \infty$. Thus $C_{0}=\{y_{\alpha}^{(n)}|0\leq\alpha\leq\kappa,\ n<\omega_{0}\}\subset C$ satisfies 
\begin{equation}\label{eq: S_{A}}
\S_{A}\subset \overline{C_{0}} 
\end{equation}
As $|C_{0}|=|\kappa\times\omega_{0}|=\kappa$ we may reindex $C_{0}=\{x_{\alpha}\}_{0\leq \alpha\leq \kappa}$. 

Since the given subgroup $\mathcal{G}^{0}$ makes $\X$ convex-transitive there exists a set of rotations 
\begin{equation}\label{eq: T0}
\{T^{(0)}_{\alpha,\beta,n}|\alpha,\beta\leq\kappa,\ n<\omega_{0}\}\subset \mathcal{G}^{0}
\end{equation}
of cardinality at most $|\kappa\times\kappa\times\omega_{0}|=\kappa$ such that for each $\alpha$ and $\beta$ it holds that 
\begin{equation}\label{eq: xbeta}
x_{\beta}\in\overline{\conv}(\{T_{\alpha,\beta,n}^{(0)}(x_{\alpha})|n<\omega_{0}\}).
\end{equation}
Note that this condition also yields that for each $a\in \overline{C_{0}}$ it holds that 
\begin{equation}\label{eq: subset}
\B_{A}\subset \overline{\conv}(C_{0})\subset\overline{\conv}(\{T^{(0)}_{\alpha,\beta,n}(a)|\alpha,\beta\leq\kappa,\ n<\omega_{0}\})
\end{equation}
by \eqref{eq: S_{A}}, \eqref{eq: T0} and \eqref{eq: xbeta}.
We may certainly assume that
\begin{equation}\label{Id}
\I_{\X}\in\{T^{(0)}_{\alpha,\beta, n}|\alpha,\beta\leq\kappa,\ n<\omega_{0}\}.
\end{equation}
 
Set 
$$A_{0}=\overline{\span}(\bigcup_{\alpha,\beta,n} T^{(0)}_{\alpha,\beta,n}(C_{0})).$$ 
Observe that $\mathrm{dens}(A_{0})=|\omega_{0}\times\kappa\times\kappa\times\omega_{0}\times\kappa|=\kappa$. 
Set 
$$B_{0}=\overline{\span}(\bigcup_{\alpha,\beta,n} (T^{(0)}_{\alpha,\beta,n})^{-1}(A_{0}))$$ 
and note that $\mathrm{dens}(B_{0})=|\omega_{0}\times\kappa\times\kappa\times\omega_{0}\times\kappa|=\kappa$. 
Fix a dense set $\{y_{\alpha}\}_{0\leq \alpha\leq \kappa}\subset \S_{B_{0}}$.
We may construct a set $C_{1}=\{z_{\gamma}\}_{0\leq \gamma\leq\kappa}\subset C$ similarly as we did $C_{0}$ above such that 
$\S_{B_{0}}\subset \overline{C_{1}}$.
There is a family of rotations 
\mbox{$\{T^{(1)}_{\alpha,\beta,n}|\alpha,\beta\leq\kappa,\ n<\omega_{0}\}\subset \mathcal{G}^{0}$} such that for each pair
$(y_{\alpha},y_{\beta})\in C_{1}\times C_{1},$ where $(\alpha,\beta)\in\kappa\times\kappa,$ we have
\[y_{\beta}\in\overline{\conv}(\{T^{(1)}_{\alpha,\beta,n}(y_{\alpha})|n<\omega_{0}\})\] 
by the convex-transitivity of $\X$ with respect to $\mathcal{G}^{0}$. Set 
\[A_{1}=\overline{\span}(\bigcup \{T^{(k)}_{\alpha,\beta,n}(C_{1})|\alpha,\beta\leq\kappa ;\ 0\leq k\leq 1;\ n<\omega_{0}\}),\]
\[B_{1}=\overline{\span}(\bigcup \{T^{(k)\ -1}_{\alpha,\beta,n}(A_{1})|\alpha,\beta\leq\kappa ;\ 0\leq k\leq 1;\ n<\omega_{0}\}).\]
Similarly as above the obtained spaces $A_{1}$ and $B_{1}$ have density characters 
$|\kappa^{3}\times\omega_{0}^{2}|=\kappa$.
Note that $A_{0}\subset B_{0}\subset A_{1}\subset B_{1}$ by \eqref{Id}.
 
We proceed the construction recursively in multiple phases as follows. 
Suppose that we have obtained the following for some $j<\omega_{0}$ and all $m\in\{0,\ldots,j\}$: 
\begin{enumerate}
\item[(i)]{a subset $C_{m+1}\subset C\ \mathrm{such\ that}\ \S_{B_{m}}\subset \overline{C_{m+1}}\ \mathrm{and}\ |C_{m+1}|=\kappa$,}
\item[(ii)]{$\{T^{(m+1)}_{\alpha,\beta,n}|\alpha,\beta\leq\kappa ,\ n<\omega_{0}\}\subset \mathcal{G}^{0}$ such that 
$$C_{m+1}\subset\overline{\conv}(\{T^{(m+1)}_{\alpha,\beta,n}(z)|\alpha,\beta\leq\kappa,\ n<\omega_{0}\})\ \mathrm{for\ each}\ z\in C_{m+1},$$}
\item[(iii)]{$A_{m+1}=\overline{\span}(\bigcup\{T^{(k)}_{\alpha,\beta,n}(C_{m+1})|\alpha,\beta\leq\kappa ;\ 0\leq k\leq m+1;\ n<\omega_{0}\})$,}
\item[(iv)]{$B_{m+1}=\overline{\span}(\bigcup\{(T^{(k)}_{\alpha,\beta,n})^{-1}(A_{m+1})|\alpha,\beta\leq\kappa ;\ 0\leq k\leq m+1;\ n<\omega_{0}\})$,}
\end{enumerate}
where the sets $C_{m+1},\ A_{m+1},\ B_{m+1}$ have density character $\kappa$. 
 
Similarly as in \eqref{eq: S_{A}} we can construct
\begin{enumerate}
\item[(v)]{$C_{j+2}\subset C\ \mathrm{such\ that}\ \S_{B_{j+1}}\subset \overline{C_{j+2}}\ \mathrm{and}\ |C_{j+2}|=\kappa$.}
\end{enumerate}
For the set $C_{j+2}$ obtained we can find similarly as in \eqref{eq: T0} and \eqref{eq: subset} a subset of rotations 
\begin{enumerate}
\item[(vi)]{$\{T^{(j+2)}_{\alpha,\beta,n}|\alpha,\beta\leq\kappa ,\ n<\omega_{0}\}\subset \mathcal{G}^{0}$ such that 
$$C_{j+2}\subset\overline{\conv}(\{T^{(j+2)}_{\alpha,\beta,n}(z)|\alpha,\beta\leq\kappa,\ n<\omega_{0}\})\ \mathrm{for\ each}\ z\in C_{j+2}.$$}
\end{enumerate}
Next we put
\begin{enumerate}
\item[(vii)]{$A_{j+2}=\overline{\span}(\bigcup\{T^{(k)}_{\alpha,\beta,n}(C_{j+2})|\alpha,\beta\leq\kappa ;\ 0\leq k\leq j+2;\ n<\omega_{0}\})$}
\item[(viii)]{$B_{j+2}=\overline{\span}(\bigcup\{(T^{(k)}_{\alpha,\beta,n})^{-1}(A_{j+2})|\alpha,\beta\leq\kappa ;\ 0\leq k\leq j+2;\ n<\omega_{0}\})$.} 
\end{enumerate}
By the induction hypothesis $\mathrm{dens}(B_{j+1})=\kappa$. 
Thus, the definition of $A_{j+2}$ and $B_{j+2}$ 
gives that $\mathrm{dens}(A_{j+2})=\mathrm{dens}(B_{j+2})=|\kappa^{3}\times\omega_{0}^{2}\times (j+1)|=\kappa$ 
similarly as above. This completes the recursive construction.

Let $\Y=\overline{\bigcup_{m<\omega_{0}}A_{m}}$. Since $A_{m}$ has density character $\kappa$ for all $m<\omega_{0}$ 
we have that $\Y$ has density character $|\kappa\times\omega_{0}|=\kappa$ also.
We aim to show that $T_{\alpha,\beta,n}^{(k)}(\Y)=\Y$ for all $\alpha,\beta \leq\kappa$ and 
$k,n<\omega_{0}$. By the construction the above family of rotations 
$\{T^{(j+1)}_{\alpha,\beta,n}\}_{\alpha,\beta,n}$ for a given $j+1<\omega_{0}$ satisfies 
\begin{equation}\label{eq: BC}
\B_{B_{j}}\subset\overline{\conv}(C_{j+1})\subset\overline{\conv}(\{T^{(j+1)}_{\alpha,\beta,n}(a)|\alpha,\beta\leq \kappa,\ n<\omega_{0}\})
\end{equation}
for each $a\in \S_{B_{j}}$ since $\S_{B_{j}}\subset \overline{C_{j+1}}$.
The facts \eqref{Id} and \eqref{eq: BC} give the chain of inclusions  
\[\S_{A}\subset\overline{C_{0}}\subset \S_{A_{0}}\subset \S_{B_{0}}\subset \overline{C_{1}}\subset \S_{A_{1}}\subset \S_{B_{1}}\subset\overline{C_{2}}\subset \S_{A_{2}}\subset \S_{B_{2}}\subset\overline{C_{3}}\subset\ldots\]
Note that from the construction of the subspaces $A_{0}, A_{1},\ldots$ and $B_{0},B_{1},\ldots$ it follows for each $m<\omega_{0}$ and $j\leq m$ that
$$T^{(j)}_{\alpha,\beta,n}(A_{m})\subset A_{m+1},\ (T^{(j)}_{\alpha,\beta,n})^{-1}(A_{m})\subset B_{m}\subset A_{m+1}.$$
Consequently
$$T^{(k)}_{\alpha,\beta,n}(\bigcup_{m<\omega_{0}}A_{m})\subset \bigcup_{m<\omega_{0}}A_{m}\ \mathrm{and}\ (T^{(k)}_{\alpha,\beta,n})^{-1}(\bigcup_{m<\omega_{0}}A_{m})\subset \bigcup_{m<\omega_{0}}A_{m}$$ 
for all $\alpha,\beta\leq\kappa $ and $n,k<\omega_{0}$. 
Hence $T^{(k)}_{\alpha,\beta,n}(\bigcup_{m<\omega_{0}}A_{m})=\bigcup_{m<\omega_{0}}A_{m}$ for all 
$\alpha,\beta\leq\kappa $ and $n,k<\omega_{0}$. 
By the uniform continuity of $T^{(k)}_{\alpha,\beta,n}$ and $(T^{(k)}_{\alpha,\beta,n})^{-1}$ we may extend
\[T^{(k)}_{\alpha,\beta,n}(\Y)=T^{(k)}_{\alpha,\beta,n}(\overline{\bigcup_{m<\omega_{0}}A_{m}})=\overline{\bigcup_{m<\omega_{0}} A_{m}}=\Y\]and similarly
\[(T^{(k)}_{\alpha,\beta,n})^{-1}(\Y)=(T^{(k)}_{\alpha,\beta,n})^{-1}(\overline{\bigcup_{m<\omega_{0}}A_{m}})=\overline{\bigcup_{m<\omega_{0}} A_{m}}=\Y\] 
for all $\alpha,\beta\leq\kappa$ and $n,k<\omega_{0}$.
 
So far we have established that 
\[\mathcal{F}=\{T^{(k)}_{\alpha,\beta,n\ |\Y}|k,n<\omega_{0};\ \alpha,\beta\leq\kappa \}\cup\{(T^{(k)}_{\alpha,\beta,n\ |\Y})^{-1}|k,n<\omega_{0};\ \alpha,\beta\leq\kappa \}\] 
is some set of rotations of $\Y$. Consider the subgroup 
\[\mathcal{G}_{\Y}^{0}=\{T_{|\Y}\in\mathcal{G}_{\Y}|\ T\in\mathcal{G}^{0}\ \mathrm{and}\ T(\Y)=\Y\}\supset \mathcal{F}.\]
 
Let $\epsilon>0$ and take two arbitrary points $x,y\in \S_{\Y}$. Then there is $m<\omega_{0}$ such that 
$\dist(x,\S_{A_{m}})<\frac{\epsilon}{2}$ and $\dist(y,\S_{A_{m}})<\frac{\epsilon}{2}$, since 
$\bigcup_{m<\omega_{0}}\S_{A_{m}}$ is dense in $\S_{\Y}$. Fix $a,b\in \S_{A_{m}}$ such that
$||x-a||<\frac{\epsilon}{2}$ and $||y-b||<\frac{\epsilon}{2}$.
Note that for $a\in A_{m}$ it holds that 
\begin{equation}\label{eq: BAT}
\B_{A_{m}}\subset \overline{\conv}(\{T(a)|T\in \mathcal{F}\}).
\end{equation} 
By applying Lemma \ref{fact} we get 
$b\in \overline{\conv}(\mathcal{G}_{\Y}^{0}(x))+\frac{\epsilon}{2}\B_{\Y}$, so that 
\[\dist(y,\overline{\conv}(\mathcal{G}_{\Y}^{0}(x)))\leq \dist(b,\overline{\conv}(\mathcal{G}_{\Y}^{0}(x)))+||y-b||\leq \frac{\epsilon}{2}+\frac{\epsilon}{2}=\epsilon.\]
Since $\epsilon>0$ was arbitrary we conclude that $y\in \overline{\conv}(\mathcal{G}_{\Y}^{0}(x))$, 
so that $\Y$ is a convex-transitive space with respect to $\mathcal{G}_{\Y}^{0}$. 
Finally, $\overline{C\cap \S_{\Y}}=\S_{\Y}$, since $\overline{\bigcup_{m<\omega_{0}} C_{m}}=\S_{\Y}$.
\end{proof}
\begin{remark}
In fact, an analogous result holds also in the almost transitive setting with a straight-forward modification of the above proof.
This result should be compared to \cite[Thm. 2.24]{BR2} courtesy of Cabello \cite{Ca4}.
\end{remark}
Recall that for each separable $\X$ there is an isometric embedding $\X\hookrightarrow \ell^{\infty}\hookrightarrow L^{\infty}(0,1)$. 
We obtain the following consequence by applying Theorem \ref{extpro} to the convex-transitive space $L^{\infty}(0,1)$. 
(Note that $\ell^{\infty}$ is not convex-transitive.)
\begin{corollary}
For each separable Banach space $\X$ there exists a separable Banach space $\Y$ such that
$\X\subset\Y\subset L^{\infty}(0,1)$ isometrically (up to suitable identifications) 
and such that $\Y$ is convex-transitive with respect to $\mathcal{G}_{\ast}=\{T_{|\Y}|\ T\in\mathcal{G}_{L^{\infty}},\ T(\Y)=\Y\}$.
\end{corollary} 

\section{Big points and Flinn elements}
Related to question ($\mathbf{Q}_{\mathrm{R}}$) one could ask if the mere existence of a big point $u\in\S_{\X}$, which is simultaneously 
a Flinn element, in fact guarantees that $\X$ is isometrically a Hilbert space. This is not the case, as the following example shows: 
\begin{example}\label{example2}
In $\ell^{1}$ the canonical unit vectors $e_{n}$ are simultaneously big points and Flinn elements, since 
$\pm e_{j}\in \mathcal{G}(e_{n})$ for $j\in\N,$ where $\B_{\ell^{1}}=\overline{\conv}(\{\pm e_{j}|j\in \N\})$ and 
$(e_{n},e_{n}^{\ast})\in \S_{\ell^{1}}\times \S_{\ell^{\infty}}$ are Flinn pairs. 
\end{example}
However, it turns out that if $\X$ meets certain (mild) geometric conditions, then this kind of characterization does actually hold.
In this section we will give such sufficient conditions.

\begin{theorem}\label{bpt}
Let $\X$ be a Banach space such that there exists a bicontractive projection $P\colon \X\rightarrow [u]$, where 
$u\in\S_{\X}$ is a big point. Assume additionally that one of the following conditions hold:
\begin{enumerate}
\item[(a)]{For every sequence $(x_{n})\subset \B_{\X}$ such that $||u+x_{n}||\rightarrow 2$ as $n\rightarrow \infty$ it holds that $P(x_{n})\rightarrow u$ as $n\rightarrow \infty$.}
\item[(b)]{The dual $\X^{\ast}$ is $\omega^{\ast}$-LUR.}
\end{enumerate}
Then $\X$ is in fact a Hilbert space. 
\end{theorem}
\begin{proof}
We write $P=f\otimes u$, where $(u,f)\in \S_{\X}\times \S_{\X^{\ast}}$ is a Flinn pair.
Let $A\subset \X$ be a $3$-dimensional subspace and fix $a\in \S_{A}$. By Proposition \ref{char} 
it suffices to show that $(a,g)$ is a Flinn pair for some suitable $g\in \S_{A^{\ast}}$.
 
We begin by constructing a separable subspace $\Y$ containing $A$ such that for the orbit $\mathcal{G}_{\X}(u)$ 
of the big point $u$ it holds that $\overline{\conv}(\mathcal{G}_{\X}(u))\cap \S_{\Y}=\S_{\Y}$.
 
Put $A_{1}=A$. Since $A_{1}$ is separable, we may approximate $\S_{A_{1}}$ by countably many elements coming from 
$\conv(\mathcal{G}_{\X}(u))$ as $u$ is a big point in $\X$.  Indeed, take a countable dense family 
$\{a_{j}\}_{j\in \N}\subset \S_{A_{1}}$ and select for each $j\in \N$ a sequence 
$\{y_{l}^{(j)}|l\in \N\}\subset \mathcal{G}_{\X}(u)$ such that $a_{j}\in \overline{\conv}(\{y_{l}^{(j)}|l\in \N\})$.  
Thus there exists a set $\{z_{n}^{(1)}|n\in \N\}=\{y_{l}^{(j)}|l,j\in\N\} \subset \mathcal{G}_{\X}(u)$ 
such that $\B_{A_{1}}$ is contained in $\overline{\conv}\{z_{n}^{(1)}|n\in \N\}$, a closed convex set containing $\S_{A_{1}}$. 
 
Put $A_{2}=\overline{\span}(\{z_{n}^{(1)}|n\in \N\})$ and note that this is again a separable subspace of $\X$. We proceed recursively:
For any $k\in \N$ the space $A_{k+1}=\overline{\span}(\{z_{n}^{(k)}|n\in \N\})$ is separable. Since $u$ is a big point in $\X$, there 
exists a countable set $\{z_{n}^{(k+1)}|n\in \N\}\subset \mathcal{G}_{\X}(u)$ such that 
$\B_{A_{k+1}}\subset \overline{\conv}(\{z_{n}^{(k+1)}|n\in \N\}$.
 
Now let $A_{1}\subset A_{2}\subset A_{3}\subset\ldots$ be the resulting sequence of subspaces. Put 
$$\Y=\overline{\bigcup_{k\in \N}A_{k}}=\overline{\span}\{z_{n}^{(k)}|n,k\in \N\}.$$
The latter equality holds because each $y\in \Y$ can be approximated by a sequence $(y_{k})\subset \Y$
such that $y_{k}\in A_{k}$ for each $k\in \N$. Since $\bigcup_{k\in \N} \S_{A_{k}}$ is dense in $\S_{\Y}$ we get
$$\B_{\Y}=\overline{\conv}\{z_{n}^{(k)}|n,k\in \N\}.$$ 
Note that $\Y$ is separable by construction.   
After this stage we will work mainly in the space $\Y$.
 
By Mazur's theorem (Thm. \ref{classmazur}) the set of smooth points of $\S_{\Y}$ relative to $\Y$ 
is dense. Let $x\in \S_{\Y}$ be a smooth point. We claim that $x$ is a Flinn element in $\Y$.
Let $g_{0}\in \S_{\Y^{\ast}}$ be the unique supporting functional for $x\in \Y$. 

Fix a sequence $(x_{m})\subset \{z_{n}^{(k)}|n,k\in \N\}$ such that $x\in \overline{\conv}\{x_{m}|m\in \N\}$.
By the construction of the set $\{z_{n}^{(k)}|n,k\in \N\}$ there is a sequence 
$(T_{m})_{m}\subset \mathcal{G}_{\X}$ such that $T_{m}(u)=x_{m}$ for $m\in \N$. 
Put $f_{m}=f\circ T_{m}^{-1}\in \X^{\ast}$ for each $m\in\N$.

There exists by Lemma \ref{lm} a sequence of finite subsets $I_{k}\subset \N$ such that  
\begin{eqnarray}
  &&\sup_{m\in I_{k}}|1-g_{0}(x_{m})|\rightarrow 0\ \mathrm{as}\ k\rightarrow\infty,\label{eq: sup}\\
  &&dist(x,\conv(\{x_{m}|m\in I_{k}\}))\rightarrow 0\ \mathrm{as}\ k\rightarrow\infty.\label{eq: dist}
\end{eqnarray}

Put $h_{m}=f_{m\ |\Y}\in \Y^{\ast}$. Let $(e_{k})\subset \B_{\Y^{\ast}}$ be an arbitrary sequence such that
$e_{k}\in \{h_{m}|m\in I_{k}\}$ for each $k\in \N$.
At this stage we proceed by applying the parallel assumptions $(a)$ and $(b)$ separately.

By applying \eqref{eq: sup} together with the fact that $g_{0}(x)=||g_{0}||=1$ we obtain 
$$\sup_{m\in I_{k}}(2-||x+x_{m}||)\leq \sup_{m\in I_{k}}(2-g_{0}(x+x_{m}))=\sup_{m\in I_{k}}(1-g_{0}(x_{m}))\rightarrow 0$$ 
as $k\rightarrow \infty$.
Note that
$$||u+T_{m}^{-1}(x)||=||T_{m}^{-1}(x_{m}+x)||=||x_{m}+x||,\ m\in\N.$$
Thus $\sup_{m\in I_{k}}(2-||u+T_{m}^{-1}(x)||)\rightarrow 0$ as $k\rightarrow \infty$. This yields that $\sup_{m\in I_{k}}(1-f(T_{m}^{-1}(x)))\rightarrow 0$ as $k\rightarrow \infty$. 
Indeed, if this were not the case then one could fix a sequence $(m_{k})_{k\in \N}$ such that $m_{k}\in I_{k}$ for each $k\in \N$ 
and $f(T_{m_{k}}^{-1}(x))\not\rightarrow 1$ as $k\rightarrow \infty$. But this is a contradiction since assumption $(a)$ yields 
that 
\[f(T_{m_{k}}^{-1}(x))=f(P(T_{m_{k}}^{-1}(x)))+f((\I-P)(T_{m_{k}}^{-1}(x)))\stackrel{k\rightarrow\infty}{\longrightarrow} f(u)=1.\]
Hence we have obtained that
$\sup_{m\in I_{k}}(1-f_{m}(x))\rightarrow 0$ as $k\rightarrow \infty$. Recall that $x$ is smooth in $\S_{\Y}$  
and note that $e_{k}(x)\rightarrow 1$. Thus we have by the \v Smulyan lemma (\ref{smulyan1}) that 
$e_{k}\stackrel{\omega^{\ast}}{\longrightarrow} g_{0}$ in $\Y^{\ast}$ as $k\rightarrow \infty$. 

Under assumption $(b)$ analogous facts can be obtained as follows. 
Let $(m_{k})\subset \N$ be a sequence such that $f_{m_{k}}\ _{|\Y}=e_{k}$.
But $f_{m_{k}}(x_{m_{k}})=f(T^{-1}_{m_{k}}(x_{m_{k}}))=f(u)=1$ by the construction for each $k$, while on the other hand 
$g_{0}(x_{m_{k}})\rightarrow 1$ as $k\rightarrow \infty$ by the construction of the sets $I_{k}$.
Under the $(b)$-assumption $\X^{\ast}$ was assumed to be $\omega^{\ast}$-LUR, and thus the condition
$$||f_{m_{k}}+g_{0}||\geq f_{m_{k}}(x_{m_{k}})+g_{0}(x_{m_{k}})\rightarrow 2\ \mathrm{as}\ k\rightarrow \infty $$
gives that $f_{m_{k}}\ _{|\Y}\stackrel{\omega^{\ast}}{\longrightarrow} g_{0}$ in $\Y^{\ast}$. 
Here the $\omega^{\ast}$-LUR property is inherited by $\Y^{\ast}$. 

We conclude that under both the assumptions $(a)$ and $(b)$ for any sequence $(e_{k})\subset \B_{\Y^{\ast}}$ such that 
$e_{k}\in \{h_{m}|m\in I_{k}\}$ for each $k\in \N$ one has that 
$e_{k}\stackrel{\omega^{\ast}}{\longrightarrow}g_{0}$ in $\Y^{\ast}$ as $k\rightarrow\infty$. 
Since $(e_{k})$ was arbitrary it must hold that
$$\sup_{m\in I_{k}}|h_{m}(y)-g_{0}(y)|\rightarrow 0$$ 
for each $y\in \Y$ as $k\rightarrow \infty$.

Let us sum up the properties of $(I_{k})$ obtained here. By the construction of the sequence $(I_{k})$ there are
$(c_{m}^{(k)})_{m}\in \S_{\ell^{1}_{+}}\cap c_{00}$ for $k\in\N$ such that
\begin{equation}  
\sum_{m\in I_{k}}c_{m}^{(k)}x_{m}\stackrel{||\cdot||}{\longrightarrow}x\ \mathrm{as}\ k\rightarrow \infty. \label{eq: ax}
\end{equation}
For each $y\in \Y$ it holds that
\begin{equation}
\sum_{m\in I_{k}}c_{m}^{(k)}|h_{m}(y)-g_{0}(y)|\leq \sup_{m\in I_{k}}|h_{m}(y)-g_{0}(y)|\sum_{m\in I_{k}}c_{m}^{(k)}\rightarrow 0\ \mathrm{as}\ k\rightarrow \infty.\label{eq: vsup}
\end{equation}
The convex combinations above in \eqref{eq: ax} are given by \eqref{eq: dist}.

Note that $(x_{m},h_{m})\in \S_{\Y}\times \S_{\Y^{\ast}}$ are Flinn pairs relative to $\Y$ for $m\in\N$.
Indeed,
\begin{eqnarray*}
||\I_{\Y}-h_{m}\otimes x_{m}||_{L(\Y)}\leq ||\I-f_{m}\otimes x_{m}||_{L(\X)}=||\I-(f\circ T_{m}^{-1})\otimes T_{m}(u)||_{L(\X)}\\
=||(T_{m}\circ T_{m}^{-1})-T_{m}\circ (f\otimes u)\circ T_{m}^{-1}||_{L(\X)}=||T_{m}\circ (\I_{\X}-f\otimes u)\circ T_{m}^{-1}||_{L(\X)}=1,
\end{eqnarray*} 
since $(u,f)$ is a Flinn pair for $\X$. Now by \eqref{eq: ax}, \eqref{eq: vsup} and Lemma \ref{lm2} we get that 
$(x,g_{0})$ is a Flinn pair relative to $\Y$.  
 
Since the smooth point $x\in \S_{\Y}$ belongs to a dense set we may apply the fact due to Kalton and Randrianantoanina 
that the set of Flinn elements is closed, see \cite[Prop. 3.1]{KR}. As a consequence we get that each point in $\S_{\Y}$ is 
actually a Flinn element relative to $\Y$. In particular $a\in \S_{A}$ is a Flinn element 
relative to $A$ and this completes the proof.
\end{proof}

It turns out in the following variation of the previous result that the $\omega$-LUR property of $\X$ provides suitable
control of the geometry, so that something more can be said. In \cite{asypri} \emph{almost isometric} versions 
(in the sense of the Banach-Mazur distance) of the existence of a Flinn element and almost transitivity
were studied. Here we apply an assumption regarding an orthogonal-like decomposition of almost isometric nature, instead 
of assuming the existence of a Flinn element. 

\begin{theorem}\label{wlurth}
Let $\X$ be a $\omega$-LUR Banach space. Assume that the following conditions hold:
\begin{enumerate}
\item[(i)]{There exists a big point $u\in \S_{\X}$.}
\item[(ii)]{For each $\epsilon>0$ there exists a $1$-codimensional, $(1+\epsilon)$-complemented subspace $\Z_{\epsilon}\subset \X$.}
\end{enumerate}
Then $\X$ is a Hilbert space.
\end{theorem}
\begin{proof}
Let $A\subset \X$ be a $3$-dimensional subspace and let $a\in \S_{A}$ be arbitrary. 
It suffices as above to show that $(a,g)$ is a Flinn pair for some suitable $g\in A^{\ast}$.

By the assumption (ii) there are pairs $(f_{k},v_{k})\in \X^{\ast}\times \X$ for $k\in\N$ such that  
$f_{k}(v_{k})=1$ and 
\begin{equation}\label{eq: redundant}
||\I-f_{k}\otimes v_{k}||\leq 1+\frac{1}{k},\ ||f_{k}\otimes v_{k}||\leq 2+\frac{1}{k}
\end{equation}
for all $k\in\N$, so that we may normalize such that $v_{k}\in \S_{\X}$ and $f_{k}\in 3\B_{\X^{\ast}}$. 
Let $e_{a}\in \S_{\X^{\ast}}$ be a support functional for $a$. 

Under the $\omega$-LUR assumption $u$ is a big point if and only if for each $x\in \S_{\X}$ there exists a sequence 
$(u_{n})\subset \mathcal{G}(u)$ such that $u_{n}\stackrel{\omega}{\rightarrow} x$ as 
$n\rightarrow \infty$. Indeed, if $h\in \S_{\X^{\ast}}$ is a support functional for $x$, then
\begin{equation}\label{eq: 1hx}
1=h(x)\in h(\overline{\conv}(G(u)))\subset \overline{\conv}(h(G(u))),
\end{equation}
so that $\sup_{z\in \mathcal{G}(u)}||x+z||\geq\sup_{z\in\mathcal{G}(u)}h(x+z)=2$. 
Thus one can select a sequence $(u_{n})\subset \mathcal{G}(u)$ such that
$||x+u_{n}||\rightarrow 2$ as $n\rightarrow\infty$. Under the $\omega$-LUR assumption on $\X$ this means that 
$u_{n}\stackrel{\omega}{\longrightarrow} x$ as $n\rightarrow\infty$. 
On the other hand, if $(y_{n})\subset \mathcal{G}(u)$ is such that $y_{n}\stackrel{\omega}{\longrightarrow}x$ as $n\rightarrow\infty$ then
Mazur's theorem yields that $x\in \overline{\conv}(\{y_{n}|n\in \N\})$. 

Since $u$ is a big point there exist families $\{S_{i}\}_{i\in \N},\{T_{i,k}\}_{(i,k)\in \N^{2}}\subset \mathcal{G}$ such that the
families $\{a_{i}\}_{i}\stackrel{\cdot}{=}(S_{i}(u))_{i},\ i\in\N,$ and 
$\{u_{i}^{k}\}_{(i,k)}\stackrel{\cdot}{=}\{T_{i,k}(u)\}_{(i,k)},\ (i,k)\in\N\times N,$ satisfy that 
$a_{i}\stackrel{\omega}{\longrightarrow} a$ and $u_{i}^{k}\stackrel{\omega}{\longrightarrow}v_{k}$ as $i\rightarrow\infty$ 
for each $k\in \N$. This means that $||a+a_{i}||\rightarrow 2$ and $||v_{k}+u_{i}^{k}||\rightarrow 2$ 
as $i\rightarrow \infty$ for each $k\in \N$. 

Let $g\in \S_{\X^{\ast}}$ be a support functional of $u$. Let $\{v_{i}^{k}\}_{(i,k)}\stackrel{\cdot}{=}\{T_{i,k}^{-1}(v_{k})\}_{(i,k)}$. 
Then 
\[||v_{i}^{k}+u||=||T_{i,k}^{-1}(v_{k}+u_{i}^{k})||=||v_{k}+u_{i}^{k}||\rightarrow 2\ \mathrm{as}\ i\rightarrow\infty\] 
and hence $g(v_{i}^{k})\rightarrow 1$ as $i\rightarrow \infty$ for all $k\in \N$ according to the $\omega$-LUR assumption.
Similarly we have for $(b_{i})_{i}\stackrel{\cdot}{=}(S_{i}^{-1}(a))_{i},\ i\in\N,$ that $g(b_{i})\rightarrow 1$ as $i\rightarrow \infty$. 
This means that $g(b_{i}+v_{i}^{k})\rightarrow 2$ as $i\rightarrow \infty$ for all $k\in\N$. 
Thus $||b_{i}+v_{i}^{k}||\rightarrow 2$ as $i\rightarrow\infty$ for each $k\in\N$. 
Put $\b_{i,k}=T_{i,k}(b_{i})$ and $\vv_{i,k}=S_{i}(v_{i}^{k})$ for all $i,k\in \N$. 
We obtain 
$$||\b_{i,k}+v_{k}||=||T_{i,k}(b_{i}+v_{i}^{k})||=||b_{i}+v_{i}^{k}||\rightarrow 2\ \mathrm{as}\ i\rightarrow \infty$$
and 
$$||a+\vv_{i,k}||=||S_{i}(b_{i}+v_{i}^{k})||=||b_{i}+v_{i}^{k}||\rightarrow 2\ \mathrm{as}\ i\rightarrow \infty.$$
Thus we get by the assumption that $\X$ is $\omega$-LUR that 
\begin{equation}\label{eq: omega}
\b_{i,k}\stackrel{\omega}{\longrightarrow} v_{k}\ \mathrm{and}\ \vv_{i,k}\stackrel{\omega}{\longrightarrow} a 
\end{equation}
as $i\rightarrow \infty$ for each $k\in \N$. Put $f_{k}^{i}=f_{k}\circ T_{i,k}\circ S_{i}^{-1}$. 
We obtain by \eqref{eq: omega} and the fact $f_{k}(v_{k})=1,\ k\in \N,$ that 
$f_{k}^{i}(a)=f_{k}(T_{i,k}(S_{i}^{-1}(a)))=f_{k}(T_{i,k}(b_{i}))=f_{k}(\b_{i,k})\rightarrow 1$ and 
$||a+\vv_{i,k}||\rightarrow 2$ as $i\rightarrow \infty$ for each $k\in \N$. Hence we may define 
\[m_{k}=\min\{j\in \N| f_{k}^{i}(a)\geq 1-\frac{1}{k}\ \mathrm{and}\ ||a+\vv_{i,k}||\geq 2-\frac{1}{k}\ \mathrm{for\ all}\ i\geq j\}\]
for each fixed $k\in\N$. Then it holds that 
\begin{equation}\label{eq: veto}
f_{k}^{m_{k}}(a)\rightarrow 1\ \mathrm{and}\ ||a+\vv_{m_{k},k}||\rightarrow 2\ \mathrm{as}\ k\rightarrow \infty.
\end{equation} 
Observe that by the $\omega$-LUR condition again
\begin{equation}\label{eq: wconv}
\vv_{m_{k},k}\stackrel{\omega}{\rightarrow}a\ \mathrm{as}\ k\rightarrow \infty. 
\end{equation}
Denote $w_{k}=\vv_{m_{k},k}$ for $k\in\N$.

Let $\Y=\overline{\span}(\{w_{k}|k\in\N\})+A$. Put $g_{k}=f^{m_{k}}_{k\ |\Y}\in 3\B_{\Y^{\ast}}$ for $k\in\N$. 
Since $\Y$ is separable we have by Alaoglu's theorem that $3\B_{\Y^{\ast}}$ is metrizable and compact in the $\omega^{\ast}$-topology. 
Thus there exists a subsequence $(k_{n})$ such that $(g_{k_{n}})$ is $\omega^{\ast}$-convergent in $3\B_{\Y^{\ast}}$. 
Denote $g_{0}=\omega^{\ast}$-$\lim_{n\rightarrow\infty}g_{k_{n}}$. Clearly $g_{0}(a)=1$, since $\lim_{k\rightarrow\infty}g_{k}(a)=1$. 

The fact that $w_{k}\stackrel{\omega}{\rightarrow}a$ as $k\rightarrow\infty$ together with Mazur's theorem yields that 
$a\in \overline{\conv}(\{w_{k_{h}}|h\geq n\})$ for each $n\in \N$. Thus, there is
$\{(c_{j}^{(l)})_{j}|l\in\N\}\subset \S_{\ell_{+}^{1}}\cap c_{00}$ such that $\sum_{j\geq l}c_{j}^{(l)}=1$ and 
the corresponding convex combinations satisfy 
\[\sum_{j}c_{j}^{(l)} w_{k_{j}}\stackrel{||\cdot||}{\longrightarrow} a\ \mathrm{and}\  \sum_{j}c_{j}^{(l)} g_{k_{j}}(a)\longrightarrow 1\quad \mathrm{as}\ l\rightarrow\infty.\]
As in the proof of Lemma \ref{lm2} we see that
$P_{l}=\sum_{j}c_{j}^{(l)} g_{k_{j}}\otimes w_{k_{j}}\stackrel{\mathrm{SOT}}{\longrightarrow}g_{0}\otimes a$ as $l\rightarrow\infty$.
Recall that  
\[g_{k}\otimes w_{k}=(f_{k}\circ T_{m_{k}}\circ S_{m_{k}}^{-1})_{|\Y}\otimes (S_{m_{k}}T^{-1}_{m_{k},k}(v_{k}))\]
and $||\I-f_{k}\otimes v_{k}||\leq 1+\frac{1}{k}$ for each $k\in \N$.
Thus we get
$$||\I-g_{0}\otimes a||=\mathrm{lim sup}_{n\rightarrow \infty} ||\I-P_{n}||=1$$
as follows. Indeed, by recalling the isometries involved here we obtain
\begin{eqnarray*} 
  &&||\I_{\X}-f_{k}^{m_{k}}\otimes \vv_{m_{k},k}||_{L(\X)}=||\I_{\X}-(f_{k}\circ T_{m_{k}}\circ S_{m_{k}}^{-1})\otimes(S_{m_{k}}T_{m_{k},k}^{-1}v_{k})||_{L(\X)}\\
&=&||S_{m_{k}}T_{m_{k},k}^{-1}-f_{k}\otimes(S_{m_{k}}T_{m_{k},k}^{-1}(v_{k}))||_{L(\X)}\\
 &=&||S_{m_{k}}T_{m_{k},k}^{-1}(\I_{\X}-f_{k}\otimes v_{k})||_{L(\X)}=||\I_{\X}-f_{k}\otimes v_{k}||_{L(\X)}\leq 1+\frac{1}{k}.
\end{eqnarray*}
By using the preceding convex combinations and the facts that $\lim_{n\rightarrow\infty}k_{n}=\infty$ and $\sum_{j\geq l}c_{j}^{(l)}=1$ 
for all $l\in\N$ we get
\begin{equation*}
\begin{array}{l}
||\I-P_{l}||_{L(\Y)}=||\I_{\Y}-\sum_{j}c_{j}^{(l)}g_{k_{j}}\otimes w_{k_{j}}||_{L(\Y)}=||\sum_{j}c_{j}^{(l)}(\I_{\Y}-g_{k_{j}}\otimes w_{k_{j}})||_{L(\Y)}\\
\leq\sum_{j}c_{j}^{(l)}||\I_{\Y}-g_{k_{j}}\otimes w_{k_{j}}||_{L(\Y)}\leq\sum_{j}c_{j}^{(l)}||\I_{\X}-f_{k_{j}}^{m_{k_{j}}}\otimes\vv_{m_{k_{j}},k_{j}}||_{L(\X)}\\
=\sum_{j}c_{j}^{(l)}||\I_{\X}-f_{k_{j}}\otimes v_{k_{j}}||_{L(\X)}\leq \sum_{j}c_{j}^{(l)}(1+\frac{1}{k_{j}})\leq 1+\frac{1}{k_{l}}\rightarrow 1\ \mathrm{as}\ n\rightarrow \infty.
\end{array}
\end{equation*}
This means that $(a,g_{0}|_{A})\in \S_{A}\times \S_{A^{\ast}}$ is a Flinn pair. Since $a\in \S_{\X}$ was arbitrary we can apply
the characterization \ref{char} to conclude that $\X$ is a Hilbert space.
\end{proof} 

In fact the argument used in the previous proof gives the following result.
\begin{proposition}
Let $\X$ be a $\omega$-LUR Banach space in which there is a big point $u\in \S_{\X}$. Then $\X$ is already convex-transitive.
\begin{proof}
We will apply suitable parts of the proof of Theorem \ref{wlurth}. For instance, here one 
\emph{does not} require the norm estimates \eqref{eq: redundant}. Hence we may substitute without loss of generality 
$v_{k}=x$ for each $k\in\N$.

Let $x,a\in\S_{\X}$ be arbitrary. We claim that $a\in\overline{\conv}(\mathcal{G}(x))$. Indeed, fix a support functional 
$f\in \S_{\X^{\ast}}$ of $x$. By following the argument of Theorem \ref{wlurth} applied to $(f,x)$, instead of the sequence 
$((f_{k},v_{k}))_{k\in\N}$, one obtains a sequence $(w_{k})_{k\in\N}\subset\mathcal{G}(x)$ such that 
$w_{k}\stackrel{\omega}{\longrightarrow} a$ as $k\rightarrow\infty$. This yields the claim by Mazur's theorem.
\end{proof}
\end{proposition}
 
\section{Interlude: Atoms and Isometric Reflections}
In the Example \ref{example2} it was observed that $e_{n}\in \ell^{1}$ are simultaneously Flinn elements and big points.
Hence $\ell^{1}$ witnesses the fact that the geometric assumptions can not be completely removed in Theorem \ref{bpt}.
Next one can ask how \emph{typical} is this particular obstruction.
 
First we will introduce a concept that to some extent generalizes the notion of \emph{atoms} defined in the Banach lattice setting
to the \emph{general} Banach space setting. Recall that for a given Banach lattice $(\X,\leq)$ a point $a\in \X\setminus \{0\}$
is an \emph{atom} if $\{v\in \X:\ 0\leq v\leq |a|\}=\{\lambda |a|:\ 0\leq \lambda \leq 1\}$. Note that in this event 
there are \emph{no} disjoint $x,y\in \{v\in\X|\ 0<v\leq |a|\}$. Recall that the points $x,y\in\X$ are disjoint, 
$x\bot y$ for short, if $|x|\wedge |y|=0$. We refer to \cite{LTII} for definitions and results concerning Banach lattices.

Let $\X$ be a Banach space. When $\Y\subset \X$ is a closed subspace, we say that a continuous linear surjective 
projection $P\colon \X\rightarrow\Y$ is an \emph{isometric reflection projection} if $(\I-2P)\in \mathcal{G}_{\X}$. 
In this event we say that $\Y$ is an \emph{isometric reflection subspace}.  

We say that $u\in \S_{\X}$ is a \emph{strong atom} if the following conditions hold:
\begin{enumerate}
\item[(1)]{$[u]$ is an isometric reflection subspace.}
\item[(2)]{For all closed subspaces $\Y\subset \X$ and all isometric reflection projections $P\colon\X\rightarrow\Y$ 
either $P(u)=u$ or $P(u)=0$.}
\end{enumerate}

Let us consider some basic examples of strong atoms. This topic will be discussed in more detail in \cite{IR}. 
Suppose that $\X$ is a $\sigma$-complete Banach lattice and $u\in\S_{\X}$ is a strong atom. 
Then $u$ is an atom also in the classical sense. This can be seen by using the existence of suitable band projections provided by the 
$\sigma$-completeness, see \cite[p.8]{LTII}.
For $\ell^{p},\ 1\leq p<\infty,\ $ the canonical unit basis vectors $e_{n}$ for any $n\in\N$ 
are strong atoms if $p\neq 2$. This is due to the characterization of rotations of $\ell^{p}$ (see e.g. \cite[2f.14]{LTI}) 
as isometric reflections are rotations as well. In a Hilbert space $\H,\ \dim(\H)\geq 2,$ there are no strong atoms $u\in\S_{\H}$. 
This is seen by fixing a suitable orthonormal basis. 

Recall that isometric reflection vectors have been used in characterizing Hilbert spaces (see \cite{SZ}).
The next theorem indicates that $\ell^{1}(\Gamma)$ is in some sense a typical obstruction for the straightforward generalization 
of ($\mathbf{Q}_{\mathrm{R}}$) to the big point setting.
\begin{theorem}
Let $\X$ be a Banach space, which has a strong atom $u\in\S_{\X}$, which is also a big point.
Then $\X=\el^{1}(\Gamma)$ isometrically for some set $\Gamma$.
\begin{proof}
Denote $\mathcal{A}=\{z\in \S_{\X}|z\ \mathrm{is\ a\ strong\ atom}\}$. Clearly the strong atom property is invariant under rotations.
Hence if $z\in\mathcal{G}(u)$ then $z\in\mathcal{A}$ and also $-z\in\mathcal{A}$.
Since $u$ is a big point we have $\overline{\conv}(\mathcal{A})=\B_{X}$. Let $(z_{k})_{k=1}^{n}\subset \mathcal{A}$ 
be a linearly independent family for some $n\in\N$. Denote the corresponding
isometric reflection projections by $P_{k}\colon \X\rightarrow [z_{k}]$. Note that since $z_{k}$ are strong atoms it follows 
that $(P_{k})_{k=1}^{n}$ are the $1$-dimensional projections associated to the basis $(z_{k})_{1}^{n}$ of 
$E\stackrel{\cdot}{=}[z_{1},\ldots ,z_{n}]$. If $z\in\mathcal{A}$ and $k\in\{1,\ldots,n\}$ are such that $P_{k}(z)\neq 0$, 
then by the definition of a strong atom we obtain that $P_{k}(z)=z$ and in fact that $z=z_{k}$. 
Put $\Gamma=\{\{\pm z\}|z\in \mathcal{A}\}$. Note that 
\begin{equation}\label{eq:zA}
\{z\in\mathcal{A}|z\neq \pm z_{k},\ 1\leq k\leq n\}\subset \bigcap_{1\leq k\leq n}\ker P_{k}.
\end{equation} 
It follows that 
\begin{equation}\label{eq: SEsub}
\S_{\X}\subset\conv\big(\conv(\{\pm z_{k}|1\leq k\leq n\})\cup \overline{\conv}(\{z\in\mathcal{A}|z\neq \pm z_{k},\ 1\leq k\leq n\})\big),
\end{equation}
where $\overline{\conv}(\{z\in\mathcal{A}|z\neq \pm z_{k},\ 1\leq k\leq n\})\subset \bigcap_{1\leq k\leq n}\ker P_{k}$ by 
\eqref{eq:zA}. 
We obtain by using \eqref{eq:zA}, a bicontractive projection $\sum_{k=1}^{n}P_{k}$ onto $E$. 
A simple argument gives that $\S_{E}\subset \conv(\{\pm z_{k}|1\leq k\leq n\})$. 

Fix $y\in \S_{E}$. Write $y=\sum_{k=1}^{n}a_{k}z_{k}\in\conv(\{\pm z_{k}|1\leq k\leq n\})$. 
Observe that above $\sum_{k=1}^{n}|a_{k}|=\sum_{k=1}^{n}||P_{k}(y)||=1$ by the triangle inequality. 
Hence $\sum_{k=1}^{n}||P_{k}(x)||=1$ for all $x\in \S_{E}$. On the other hand, suppose
that $x\in E$ is such that $\sum_{k=1}^{n}||P_{k}(x)||=1$. Then $\frac{x}{||x||}\in\S_{E}$ and thus 
\[\frac{\sum_{k=1}^{n}||P_{k}(x)||}{||x||}=\sum_{k=1}^{n}\left|\left|P_{k}\left(\frac{x}{||x||}\right)\right|\right|=1,\] 
where $\sum_{k=1}^{n}||P_{k}(x)||=1$, so that necessarily $||x||=1$. Hence we obtain 
\[||x||=\sum_{k=1}^{n}||P_{k}(x)||\ \mathrm{for\ all}\ x\in E.\] 
Consequently $E=\el^{1}_{n}$ isometrically. Moreover, since $n$ and $(z_{k})_{1}^{n}\subset \mathcal{A}$
were arbitrary we obtain that for any linearly independent subset $\{y_{k}\}_{k=1}^{m}\subset\mathcal{A}$, where $m\in\N$
it holds that $\span(\{y_{k}\}_{k=1}^{m})=\el^{1}_{m}$ isometrically.

Let us focus on the relationship between the dense subspaces $\span(\mathcal{A})\subset \X$ and 
$\span(\{e_{\gamma}\}_{\gamma\in\Gamma})\subset \ell^{1}(\Gamma)$. Consider a linear map
$F\colon \span(\mathcal{A})\stackrel{\mathrm{onto}}{\longrightarrow} \span(\{e_{\gamma}\}_{\gamma\in\Gamma})$ given by
$F([z])=[e_{\gamma}]$ for $z\in\mathcal{A}$ and $\gamma=\{\pm z\}$. We have above actually that $F$ is an isometry. 
There exists a continuous extension $\tilde{F}\colon \overline{\span}(\mathcal{A})\rightarrow \ell^{1}(\Gamma)$
of $F$, which is also a linear onto isometry. This completes the claim that $\X$ and $\ell^{1}(\Gamma)$ are mutually isometric.
\end{proof}
\end{theorem}

Consider Banach lattices $\X$ and $\Y$. Recall that an operator $T\colon \X\rightarrow\Y$ is said to be \emph{disjointness-preserving}
if $T(x)$ and $T(y)$ are disjoint whenever $x,y\in\X$ are disjoint. Suppose that a given positive measure $\mu$ has an atom and for some 
$1\leq p\leq \infty,\ p\neq 2,$ the space $L^{p}(\mu)$ is convex-transitive. Then $L^{p}(\mu)$ is actually $1$-dimensional, see 
\cite[Cor.3.5]{BR1}. We generalize this fact in the following result.
\begin{theorem}
Let $\X$ be a convex-transitive Banach lattice such that each rotation $T\in\mathcal{G}_{\X}$ is disjointness-preserving.
If there exists an atom $x\in \S_{\X}$, then $\X$ is in fact $1$-dimensional.
\end{theorem}
\begin{proof}
If there is $x\in\S_{\X}$, then $\dim(\X)\geq 1$. Assume to the contrary that $\dim(\X)\geq 2$.
Since each rotation $T\in\mathcal{G}_{\X}$ is disjointness-preserving we obtain that $\mathcal{G}_{\X}(x)$
consists of atoms. By applying the assumption that $x$ is a big point we have that $\X$ is purely atomic, that is,
there exist pairwise disjoint atoms $\{x_{\gamma}\}_{\gamma\in\Gamma}\subset \mathcal{G}_{\X}(x)$ such that 
$\overline{\span}(\{x_{\gamma}|\gamma\in\Gamma\})=\X$. Since $\dim(\X)\geq 2$ we have $|\Gamma|\geq 2$. Pick disjoint 
$x_{\gamma_{1}},x_{\gamma_{2}}\in \mathcal{G}_{\X}(x)$, where $\gamma_{1},\gamma_{2}\in\Gamma$. 
Put $y_{1}=\frac{|x_{\gamma_{1}}|}{||x_{\gamma_{1}}+x_{\gamma_{2}}||}$ and $y_{2}=\frac{|x_{\gamma_{2}}|}{||x_{\gamma_{1}}+x_{\gamma_{2}}||}$. Observe that $y_{1},y_{2}$ are positive, disjoint and $||y_{1}+y_{2}||=1$. 
Put $y=y_{1}+y_{2}$. Since $y_{1}$ and $y_{2}$ are positive we get $||y_{j}||\leq ||y||$ for $j=1,2$. 

Let us check for convenience that there exist band projections $P_{y_{1}},P_{y_{2}}$ and $P_{x}$ onto $[y_{1}], [y_{2}]$ and $[x]$, 
respectively. Indeed, by the norm continuity of the operation $v\mapsto |x|\wedge |v|$ (see e.g. \cite[p.1]{LTII}) we obtain that 
$x\notin x^{\bot}=\overline{\span}(\{v\in \X|v\bot x\})$. Since $\X$ is purely atomic, 
it holds that $\X=[x]\oplus x^{\bot}$. Note that $||x-z||\geq ||x||=1$ whenever $z\in x^{\bot}$.
Thus by the Hahn-Banach theorem there is $f\in \X^{\ast}$ such that $\ker(f)=x^{\bot}$ and $f(|x|)=||f||=||x||=1$. 
Hence $P_{x}=f\otimes |x|$ is the claimed projection. Similarly we obtain band projections $P_{y_{1}}$ and $P_{y_{2}}$. 
Clearly $||P_{x}||=||P_{y_{1}}||=||P_{y_{2}}||=1$.

Since the rotations $T\in\mathcal{G}_{\X}$ preserve atoms, we obtain that $P_{y_{1}}(v)=0$ or 
$(\I-P_{y_{1}})(v)=0$ holds for each $v\in \mathcal{G}_{\X}(x)$. Since $y_{1}\bot y_{2}$, we obtain in particular that 
$\min(||P_{y_{1}}(v)||,||P_{y_{2}}(v)||)=0$ for all $v\in \mathcal{G}_{\X}(x)$. On the other hand, 
$\max(||P_{y_{1}}(v)||,||P_{y_{2}}(v)||)\leq 1$ for $v\in \mathcal{G}_{\X}(x)$ as $||P_{y_{1}}||=||P_{y_{2}}||=||v||=1$. 
Consequently, it follows for the convex combinations that 
\begin{equation}\label{eq: Py1}
||P_{y_{1}}(z)||+||P_{y_{2}}(z)||\leq 1,\ \quad \mathrm{for\ all}\ z\in \overline{\conv}(\mathcal{G}_{\X}(x)). 
\end{equation}
Next we will apply the following facts: $y_{1}=P_{y_{1}}(y_{1})\neq 0\neq P_{y_{2}}(y_{2})=y_{2}$, 
$y=y_{1}+y_{2}\in \overline{\conv}(\mathcal{G}(x))$ as $x$ is a big point and $a+b=\min(a,b)+\max(a,b)$ for $a,b\in\R$. 
By \eqref{eq: Py1} applied for $z=y$ we get
\begin{equation*}
\begin{array}{rl}
&\max(||y_{1}||,||y_{2}||)=\max(||P_{y_{1}}(y_{1}+y_{2})||,||P_{y_{2}}(y_{1}+y_{2})||)\\
=&\max(||P_{y_{1}}(y_{1})||,||P_{y_{2}}(y_{2})||)\leq 1-\min(||P_{y_{1}}(y_{1})||,||P_{y_{2}}(y_{2})||)<1. 
\end{array}
\end{equation*}
Recall that $y$ is a big point, since $\X$ is convex-transitive.
Since all $T\in\mathcal{G}_{\X}$ are disjointness-preserving and $y_{1}\bot y_{2}$ we get 
$|T(y_{1})|\bot |T(y_{2})|$. As $x$ is an atom we get that $|T(y_{1})|\wedge |x|=0$ or $|T(y_{2})|\wedge |x|=0$. Hence 
\[||P_{x}(T(y))||=\max(||P_{x}(T(y_{1}))||,||P_{x}(T(y_{2}))||)\leq\max(||y_{1}||,||y_{2}||)\] 
for all $T(y)\in\mathcal{G}_{\X}(y)$. Thus $||P_{x}(z)||\leq \max(||y_{1}||,||y_{2}||)$ for all 
$z\in\overline{\conv}(\mathcal{G}_{\X}(y))$. This contradicts the facts that $||P_{x}(x)||=||x||=1$ and 
$x\in\overline{\conv}(\mathcal{G}(y))$.
\end{proof}

\section{The Convex-transitive case}
The following two theorems are the main results of this article. It turns out below that rather weak geometric conditions
guarantee that a convex-transitive space, which in addition admits a $1$-dimensional bicontractive projection, is isometrically a Hilbert
space. 
 
The duality mapping is denoted by $J\colon \S_{\X}\rightarrow \mathcal{P}(\S_{\X^{\ast}})$. This is the set-valued mapping defined by  
\[J(x)=\{f\in \S_{\X^{\ast}}|f(x)=1\},\quad \mathrm{for}\ x\in \S_{\X}.\] 
Recall that for topological spaces $X,Y$ and a point $x\in X$ the set-valued map $f\colon X\rightarrow \mathcal{P}(Y)$
is called \emph{upper semi-continuous} (u.s.c.) at $x$ if for each open set $V\subset Y$, which contains $f(x)$, there exists
an open neighbourhood $U\subset X$ of $x$ such that $f(U)\subset V$.

\begin{theorem}\label{ct}
Let $\X$ be a convex-transitive Banach space with a bicontractive projection $P\colon \X\rightarrow [u]$, where
$u\in\S_{\X}$. Consider $u\in \S_{\X}\subset \S_{\X^{\ast\ast}}$ and assume that the following conditions hold:
\begin{enumerate}
\item[(i)]{$\wexp(\B_{\X})\subset \S_{\X}$ is dense.}
\item[(ii)]{The set-valued map $u\circ J\colon (\S_{\X},\omega)\rightarrow \mathcal{P}([-1,1])$ is u.s.c. at $u$.}
\end{enumerate} 
Then $\X$ is a Hilbert space.
\end{theorem}
\begin{proof}
Write $P=f\otimes u$, where $(u,f)\in \S_{\X}\times \S_{\X^{\ast}}$ is a Flinn pair.
As in Theorem \ref{bpt} it suffices to show that for an arbitrary 3-dimensional subspace $A\subset \X$ any 
$x\in \S_{A}$ is a Flinn element relative to $A$. Let $A$ be such a subspace.
 
By applying Theorem \ref{extpro} to $A+[u]$ and $C=\wexp(\B_{\X})$ we obtain a separable convex-transitive subspace 
$\Y\subset \X$ such that $A+[u]\subset \Y$ and $\overline{\wexp(\B_{\X})\cap \S_{\Y}}=\S_{\Y}$. 
Note that any point $y\in \wexp(\B_{\X})\cap \S_{\Y}$ is a weakly exposed point of $\B_{\Y}$ in $\Y$.

One can verify by using the prior equation \eqref{eq: 1hx} that for a weakly exposed point 
$v\in \overline{\conv}(\mathcal{G}_{\X}(y))$, where $y\in\S_{\X}$, there is a sequence
$(T_{i})_{i}\subset\mathcal{G}_{\X}$ such that $T_{i}(y)\stackrel{\omega}{\longrightarrow} v$ as $i\rightarrow\infty$.   

There exists by assumption (ii) $\omega$-open neighbourhoods $V_{l}\subset \S_{\Y}$ of $u$ for $l\in\N$ 
such that $\{g(u)| g\in J(v)\}\subset (1-l^{-1},1]$ for each $v\in V_{l}$. Moreover, observe that
$W_{l}\stackrel{\cdot}{=}\{v\in V_{l}|f(v)>1-l^{-1}\}\subset V_{l}$ for $l\in\N$ are $\omega$-open neighbourhoods of $u$.  

Since $\Y$ is separable, Mazur's theorem (Thm. \ref{classmazur}) states that the set of smooth points is dense in $\S_{\Y}$. 
Let $b\in \S_{\Y}$ be an arbitrary smooth point. Let 
$(z_{n})\subset \S_{\Y}$ be a sequence of weakly exposed points relative to $\Y$ such that $||z_{n}-b||\rightarrow 0$ as $n\rightarrow\infty$. 
Let $(h_{n})\subset \S_{\Y^{\ast}}$ be the corresponding weakly exposing functionals. 
Since the weakly exposed points are norm dense in $\S_{\Y}$, we may a fortiori pick a sequence $(y_{i})_{i}\subset \wexp(\B_{\X})$
such that $y_{i}\in W_{i}$ for all $i\in\N$ and $y_{i}\rightarrow u$ as $i\rightarrow\infty$. 

As $Y$ is convex-transitive, there is a family 
$\{T_{n,i,k}\}_{(n,i,k)}\subset \mathcal{G}_{\Y}$ such that 
\[\omega-\lim_{k\rightarrow\infty}T_{n,i,k}(z_{n})=y_{i}\quad \mathrm{for\ all}\ n,i\in\N.\]
Since $y_{i}\rightarrow u$ as $i\rightarrow\infty$ we get
\[\lim_{i\rightarrow\infty}\omega-\lim_{k\rightarrow\infty}T_{n,i,k}(z_{n})=u\quad \mathrm{for\ all}\ n\in\N.\]
Define $A_{n,l}=\{(i,k)\in\N\times\N|T_{n,i,k}(z_{n})\in W_{l}\}$ for each fixed $n,l\in\N$. Then
\[u\in \overline{\{T_{n,i,k}(z_{n})|\ (i,k)\in A_{n,l}\}}^{\omega}\quad \mathrm{for\ all}\ n,l\in\N.\]
Put $h_{n,i,k}=h_{n}\circ T_{n,i,k}^{-1}$ for $n,i,k\in \N$.

Observe that for all $n,l\in\N$ and $(i,k)\in A_{n,l}$ we get by the definition of $W_{l}$ that
\begin{equation}\label{eq: 5.2}
f_{|\Y}(T_{n,i,k}(z_{n}))>1-\frac{1}{l},
\end{equation}
\begin{equation}\label{eq: 5.3}
\inf u\circ J(T_{n,i,k}(z_{n}))>1-\frac{1}{l}
\end{equation}
for $l\in\N$. Let $g_{n,i,k}=f_{|\Y}\circ T_{n,i,k}$ for $n,i,k\in\N$. Equation \eqref{eq: 5.2} gives that
\begin{equation}\label{eq: 5.4}
\inf_{(i,k)\in A_{n,l}}g_{n,i,k}(z_{n})\rightarrow 1\ \mathrm{as}\ l\rightarrow\infty\ \mathrm{for\ each}\ n\in\N.
\end{equation}
Note that $(h_{n}\circ T_{n,i,k}^{-1})T_{n,i,k}(z_{n})=1$, so that $h_{n}\circ T_{n,i,k}^{-1}\in J(T_{n,i,k}(z_{n}))$ for all $n,i,k\in\N$. By \eqref{eq: 5.3} we obtain
\[\inf_{(i,k)\in A_{n,l}}h_{n}(T^{-1}_{n,i,k}(u))=\inf_{(i,k)\in A_{n,l}}h_{n,i,k}(u)\rightarrow 1\ \mathrm{as}\ l\rightarrow\infty\
\mathrm{for\ each}\ n\in\N.\]

Fix a family $\{\boldsymbol{a}_{n,l}\}_{(n,l)\in\N^{2}}\subset\N^{3}$ such that $\boldsymbol{a}_{n,l}\in \{n\}\times A_{n,l}$ 
for all $n,l\in\N$. Since $h_{n}$ are weakly exposing functionals for $z_{n}$ respectively, we get that 
$T^{-1}_{\boldsymbol{a}_{n,l}}(u)\stackrel{\omega}{\longrightarrow}z_{n}$ as $l\rightarrow\infty$ for each $n\in\N$.
Mazur's theorem yields that 
\begin{equation}\label{eq: 5.5}
z_{n}\in\overline{\conv}(\{T_{\boldsymbol{a}_{n,l}}^{-1}(u)|l\geq j\})\quad \mathrm{for\ each}\ n,j\in\N.
\end{equation}
Put $B_{n,m}=\bigcup_{l\geq m}\{n\}\times A_{n,l}$ for all $n,m\in\N$.

Recall that $||g_{n,i,k}||=1$ for $n,i,k\in \N$. Hence by \eqref{eq: 5.4} we obtain
\begin{equation}
\lim_{m\rightarrow\infty}\inf_{\boldsymbol{b}\in B_{n,m}}g_{\boldsymbol{b}}(b)\geq \lim_{m\rightarrow\infty}\inf_{\boldsymbol{b}\in B_{n,m}}(g_{\boldsymbol{b}}(z_{n})-||z_{n}-b||)=1-||z_{n}-b||
\end{equation}
for each $n\in\N$. Hence one can fix a sequence $(k_{n})_{n\in\N}\subset\N$ such that 
\begin{equation}\label{eq: uusi}
\inf_{\boldsymbol{b}\in B_{n,m}}g_{\boldsymbol{b}}(b)\geq 1-||z_{n}-b||-2^{-n}\quad \mathrm{for\ all}\ m\geq k_{n},\ n\in\N.
\end{equation} 

Put $C_{n}=\bigcup_{l\geq \max(n,k_{n})}B_{n,l}$ for each $n\in\N$. 
By using \eqref{eq: 5.5} we may fix convex weights $\{c_{j,i,k}^{(n)}\}_{j,i,k}\in \S_{\ell^{1}_{+}(C_{n})}\cap c_{00}(C_{n})$ 
such that $||\sum_{C_{n}}c^{(n)}_{j,i,k}T^{-1}_{j,i,k}(u)-z_{n}||<2^{-n}$ for each $n\in\N$.
We obtain by using $\lim_{n\rightarrow\infty}z_{n}=b$, \eqref{eq: uusi} and the definition of $(C_{n})$ that
\begin{equation}\label{eq: 5.8}
\lim_{n\rightarrow\infty} \sum_{(j,i,k)\in C_{n}} c_{j,i,k}^{(n)}g_{j,i,k}(b)=1.
\end{equation}
The selection of the convex weights $c_{j,i,k}^{(n)}$ and $\lim_{n\rightarrow\infty}z_{n}=b$ gives that
\begin{equation}\label{eq: 5.9}
\lim_{n\rightarrow\infty}\sum_{(j,i,k)\in C_{n}}c_{j,i,k}^{(n)}T^{-1}_{j,i,k}(u)=b.
\end{equation}

Properties \eqref{eq: 5.8} and \eqref{eq: 5.9} allow us to apply Lemma \ref{lm2} to the sequence of
finite rank operators given by
\[P_{n}=\sum_{(j,i,k)\in C_{n}} c^{(n)}_{j,i,k}g_{j,i,k}\otimes T_{j,i,k}(u)\quad \mathrm{for}\ n\in\N.\]
Since $b$ is a smooth point, we get from Lemma \ref{lm2} that $P_{n}\rightarrow b^{\ast}\otimes b$ as $n\rightarrow\infty$
in the weak operator topology in $L(\Y)$, where $b^{\ast}\in \S_{\Y}$ is the unique support functional of $b$.
Moreover, since clearly $(T^{-1}_{j,i,k}(u),f\circ T_{j,i,k})$ are Flinn pairs for each $(j,i,k)\in \N^{3}$, 
we obtain further that $(b,b^{\ast})$ is a Flinn pair. 
Finally, by applying the Kalton-Randrianantoanina result as in the proof of Theorem \ref{bpt} we complete the proof.
\end{proof}
One can obtain the following example by an application of Theorem \ref{ct}. 
\begin{example}
Let $\X$ be a convex-transitive Banach space with a bicontractive projection $f\otimes u\colon \X\rightarrow [u]$,
where $(f,u)\in \S_{\X^{\ast}}\times \S_{\X}$. If $\X^{\ast}$ is smooth and $f$ is a $\omega^{\ast}$-LUR point, 
then $\X$ is isometrically a Hilbert space.
\end{example}
Indeed, each point $x\in\S_{\X}$ is weakly exposed, since $\X^{\ast}$ is smooth, (see e.g. \cite{ZZ}). 
On the other hand, by applying the fact that $f$ is a $\omega^{\ast}$-LUR element, one can see that $Ju=\{f\}$ and 
$J\colon (\S_{\X},\omega)\rightarrow \mathcal{P}(\S_{\X^{\ast}},\omega^{\ast})$ is u.s.c. at $u$.
Note that $u\colon (\S_{\X^{\ast}},\omega^{\ast})\rightarrow [-1,1],\ g\mapsto g(u)$ is continuous. 
It follows that $u\circ J$ is u.s.c. at $u$. 

\begin{theorem}\label{th:wlurpoint}
Let $\X$ be a convex-transitive Banach space which admits a $1$-dimensional bicontractive projection $P\colon \X\rightarrow L$.
Assume that there is a $\omega$-LUR point $x\in \S_{\X}$. Then $\X$ is a Hilbert space.
\end{theorem}

Let us make a brief philosophical remark about the techniques applied in the proof of this main result.
The key ingredient below is a kind of uniform control of the weak topology, $\phi$. This is obtained by using the interplay between 
the weak norm geometry and the rotations.

\begin{proof}
Write $P=f\otimes u\colon \X\rightarrow [u]$, where $(u,f)$ is a Flinn pair. As in the previous proofs it suffices to show that 
for any $3$-dimensional subspace $A\subset \X$ all points $a\in \S_{A}$
are Flinn elements relative to $A$. By Theorem \ref{extpro} there is a separable convex-transitive subspace 
$\Y\subset \X$ containing $A+[x,u]$. Mazur's theorem (see Thm. \ref{classmazur}) yields that the set of smooth points in
$\S_{\Y}$ relative to $\Y$ is dense. Let $b\in \S_{\Y}$ be any smooth point relative to $\Y$.     
Hence there is unique $b^{\ast}\in \S_{\Y^{\ast}}$ such that $b^{\ast}(b)=1$.

By using the convex-transitivity of $\Y$ and condition $($iii$)$ of Lemma \ref{lm} applied to 
$u\in \overline{\conv}(\mathcal{G}_{\Y}(x)),$ where $f(u)=1$, there is a sequence $(S_{k})\subset\mathcal{G}_{\Y}$ such that 
\begin{eqnarray}\label{eqa: sku}
  &&f(S_{k}(x))\rightarrow 1\ \mathrm{as}\ k\rightarrow \infty,\\
  &&u\in\overline{\conv}(\ \{S_{k}(x)|k\geq l\}\cup\{S_{k}(x)|f(S_{k}(x))=1\}\ )\ \mathrm{for\ each}\ l\in \N.
\end{eqnarray}   

Also, by the same argument applied to $b\in\overline{\conv}(\mathcal{G}_{Y}(x))$, where $b^{\ast}(b)=1$ we can find a sequence 
$(T_{n})\subset\mathcal{G}_{\Y}$ such that 
\begin{eqnarray}
  &&b^{\ast}(T_{n}(x))\rightarrow 1\ \mathrm{as}\ n\rightarrow \infty\label{eq: btn1}\\
  &&b\in\overline{\conv}(\{T_{n}(x)|n\geq l\}\cup\{T_{n}(x)|b^{\ast}(T_{n}(x))=1\})\ \mathrm{for\ each}\ l\in \N.
\end{eqnarray}
Let $h\in \S_{\Y^{\ast}}$ be such that $h(x)=1$. Since $x$ is an $\omega$-LUR point, it follows that $h$ weakly exposes $x$.
We find weakly exposing functionals for $x_{k}\stackrel{\cdot}{=}S_{k}(x)$ and $z_{n}\stackrel{\cdot}{=}T_{n}(x)$ by putting 
$h_{k}=h\circ S_{k}^{-1}$ and $g_{n}=h\circ T_{n}^{-1}$, respectively, for each $n,k\in \N$.

Since $b^{\ast}(z_{n})=b^{\ast}(T_{n}(x))\rightarrow 1$ as $n\rightarrow\infty$ and 
$f(x_{k})=f(S_{k}(x))\rightarrow 1$ as $k\rightarrow\infty$, where 
$||b^{\ast}||=||f||=1$, we get
\[b^{\ast}(b+z_{n})\leq ||b+z_{n}||=||T_{n}^{-1}(b)+T_{n}^{-1}(z_{n})||=||T_{n}^{-1}(b)+x||\rightarrow 2\ \mathrm{as}\ n\rightarrow \infty\]
and
\[f(u+S_{k}(x))\leq ||u+S_{k}(x)||=||S_{k}^{-1}(u)+x||\rightarrow 2\ \mathrm{as}\ k\rightarrow \infty,\]
where we applied the convergence of the lower bounds.   
Since $x$ is a $\omega$-LUR point, we obtain that $T_{n}^{-1}(b)\stackrel{\omega}{\longrightarrow} x$ as $n\rightarrow\infty$ and 
$S_{k}^{-1}(u)\stackrel{\omega}{\longrightarrow}x$ as $k\rightarrow \infty$. Thus
\begin{equation}\label{eq: hST}
h(T_{n}^{-1}(b))\rightarrow 1\ \mathrm{as}\ n\rightarrow\infty\ \mathrm{and}\ h(S_{k}^{-1}(u))\rightarrow 1\ \mathrm{as}\ n\rightarrow\infty\mathrm{.}
\end{equation}

Put 
\[\phi_{t}=\inf\{h(T_{n}^{-1}(b))|n\in \N\ \mathrm{satisfies}\ b^{\ast}(z_{n})\geq t\}\ \mathrm{for}\ t\in (0,1).\]
Clearly $\phi_{t}$ is non-decreasing with respect to $t$. If $\lim_{t\rightarrow 1^{-}}\phi_{t}=a<1$,
then $a\in [0,1]$ is a cluster point of $(h(T_{n}^{-1}(b)))_{n}$ by the definition of $\phi_{t}$. 
Thus, by the sequential compactness of $[0,1]$ there is a subsequence $(h(T_{n_{j}}^{-1}(b)))_{j}$ such that 
$h(T_{n_{j}}^{-1}(b))\rightarrow a$ as $j\rightarrow\infty$. This contradicts \eqref{eq: hST}.
Thus $\lim_{t\rightarrow 1^{-}}\phi_{t}=1$.

Since $S_{k}^{-1}(u)\stackrel{\omega}{\longrightarrow} x$ as $k\rightarrow \infty$, Mazur's classical theorem yields that 
$x\in \overline{\conv}(\{S_{k}^{-1}(u)|k\geq l\})$ for each $l\in \N$. Recall that 
$f(x_{k})\rightarrow 1$ as $k\rightarrow\infty$ by \eqref{eqa: sku}. 
Hence there is a finite subset $L_{m}\subset \N$ for each $m\in \N$
such that the following conditions hold:
\begin{eqnarray}\label{eqa: dfx}
\dist(x,\conv(\{S_{k}^{-1}(u)|k\in L_{m}\}))<\frac{1}{m}\\
f(x_{k})=f_{|\Y}\circ S_{k}(x)>1-\frac{1}{m}\ \mathrm{for}\ k\in L_{m}
\end{eqnarray}
We claim that for each $m\in \N$ there is $0<\beta_{m}<1$ such that the respective slices satisfy
\[\mathcal{S}_{h,\beta_{m}}=\{y\in \B_{\Y}|h(y)>\beta_{m}\}\subset \bigcap_{k\in L_{m}}\mathcal{S}_{f\circ S_{k},1-\frac{1}{m}}=\bigcap_{k\in L_{m}}\{y\in \B_{\Y}|f\circ S_{k}(y)>1-\frac{1}{m}\}.\]
Indeed, assume to the contrary that $\{y\in \B_{\Y}|h(y)>\beta\}$ intersects the complement of 
$\bigcap_{k\in L_{m}}\{y\in \B_{\Y}|f\circ S_{k}(y)>1-\frac{1}{m}\}$ for each choice of $0<\beta<1$. 
As $h$ is a weakly exposing functional of $x$ one can choose a sequence $(y_{n})\subset \B_{\Y}$ such that 
$y_{n}\stackrel{\omega}{\longrightarrow} x$ as $n\rightarrow\infty$ and 
$y_{n}\notin\bigcap_{k\in L_{m}}\{y\in \B_{\Y}|f\circ S_{k}(y)>1-\frac{1}{m}\}$ for each $n\in \N$. 
On the other hand, $\bigcap_{k\in L_{m}}\{y\in \B_{\Y}|f\circ S_{k}(y)>1-\frac{1}{m}\}$ is a $\omega$-open neighbourhood of $x$
as $L_{m}$ is finite, so that this provides a contradiction. Without loss of generality we may above assume that $\beta_{m}\rightarrow 1$ as $m\rightarrow\infty$ (though this is necessarily the case anyway). 

\textit{Claim:} For each $l\in \N$ there is a finite subset $K_{l}\subset \N$ such that the following conditions hold:
\begin{eqnarray}
\dist(b,\conv(\{T_{n}(x)|n\in K_{l}\}))<\frac{1}{l}\label{eq: bkl}\\
h(T_{n}^{-1}(b))>1-\frac{1}{l}\ \mathrm{for}\ n\in K_{l}
\end{eqnarray}   
Indeed, since $\lim_{t\rightarrow 1^{-}}\phi_{t}=1$ one may select a sequence $(t_{l})\subset (0,1)$ such that
$\phi_{t_{l}}\geq 1-\frac{1}{l}$ for each $l\in\N$. By \eqref{eq: btn1} we obtain  
$$b\in\overline{\conv}(\{T_{n}(x)|n\in\N\ \mathrm{satisfies}\ b^{\ast}(T_{n}(x))\geq t_{l}\})$$
for each $l\in\N$. Thus we may select for each $l\in\N$ a finite set 
$K_{l}\subset \{n\in\N| b^{\ast}(T_{n}(x))\geq t_{l}\}$ such that $\dist(b,\conv(\{T_{n}(x)|n\in K_{l}\}))<\frac{1}{l}$.
By the definitions of $K_{l},\phi_{t}$ and $(t_{l})$ it holds that 
\[\inf\ \{h(T_{n}^{-1}(b))|n\in K_{l}\}\geq \inf\ \{h(T_{n}^{-1}(b))|n\ \mathrm{satisfies}\ b^{\ast}(T_{n}(x))\geq t_{l}\}=\phi_{t_{l}}\geq 1-\frac{1}{l}\]
for each $l\in\N$. The claim is complete. 

Since $\beta_{m}<1$ for $m\in \N$, we may fix $l_{m}\in \N$ such that $\beta_{m}\leq 1-\frac{1}{l_{m}}$. 
From the definition of the sequences $(l_{m})$, $(K_{l})$ and $(\beta_{m})$ we have 
$$\{T_{n}^{-1}(b)|n\in K_{l_{m}}\}\subset \mathcal{S}_{h,1-\frac{1}{l_{m}}}\subset \mathcal{S}_{h,\beta_{m}}\subset\bigcap_{k\in L_{m}}\mathcal{S}_{f\circ S_{k},1-\frac{1}{m}}\subset \B_{\Y}.$$
With these notations we have that for all $m\in \N$
$$(f\circ S_{k}\circ T_{n}^{-1})(b)>1-\frac{1}{m}\ \mathrm{for\ all}\ k\in L_{m}\ \mathrm{and}\ n\in K_{l_{m}}.$$
Observe that for each $m\in \N$ it holds that
\begin{equation}\label{eq: conclude}
g(b)>1-\frac{1}{m}\ \mathrm{for}\ g\in \conv(\{f_{|\Y}\circ S_{k}\circ T_{n}^{-1}|k\in L_{m};\ n\in K_{l_{m}}\}).
\end{equation}
Note further that 
\[\bigcup_{n\in K_{l_{m}}}T_{n}(\conv(\{S_{k}^{-1}(u)|k\in L_{m}\}))\subset \conv(\{T_{n}\circ S_{k}^{-1}(u)|(n,k)\in K_{l_{m}}\times L_{m}\}),\] 
so that
\begin{equation}\label{eq: csub}
\begin{array}{rl}
&\conv \left(\bigcup_{n\in K_{l_{m}}}T_{n}(\conv(\{S_{k}^{-1}(u)|k\in L_{m}\}))\right)\\
\subset&\conv(\{T_{n}\circ S_{k}^{-1}(u)|(n,k)\in K_{l_{m}}\times L_{m}\}).
\end{array}
\end{equation}
Since $T_{n}$ and $S_{k}$ are isometries for each $n,k\in\N$, we obtain by applying \eqref{eq: csub} and Lemma \ref{fact} that 
\begin{equation}\label{eqa: dc}
\begin{array}{l}
\dist(b,\conv(\{T_{n}\circ S_{k}^{-1}(u)|n\in K_{l_{m}};\ k\in L_{m}\}))\\
\leq\dist(b,\conv(\bigcup_{n\in K_{l_{m}}}T_{n}(\conv(\{S_{k}^{-1}(u)|k\in L_{m}\}))))\\
\leq \dist(b,\conv(\{T_{n}(x)|n\in K_{l_{m}}\}))+\dist(x,\conv(\{S_{k}^{-1}(u)|k\in L_{m}\}))\\
\leq \frac{1}{l_{m}}+\frac{1}{m}\longrightarrow 0\ \mathrm{as}\ m\rightarrow\infty,
\end{array}
\end{equation}
where the last inequality holds by \eqref{eq: bkl} and \eqref{eqa: dfx}. According to \eqref{eqa: dc} we may fix convex weights 
$\{c_{n,k}^{(m)}\}_{n,k}\in \S_{\ell^{1}_{+}(K_{l_{m}}\times L_{m})}\cap c_{00}(K_{l_{m}}\times L_{m})$ for each $m\in\N$ such that 
\[\sum_{K_{l_{m}}\times L_{m}}c_{n,k}^{(m)}T_{n}\circ S_{k}^{-1}(u)\stackrel{||\cdot||}{\longrightarrow}b\ \mathrm{as}\ m\rightarrow \infty.\]
By \eqref{eq: conclude} we obtain 
\[\sum_{K_{l_{m}}\times L_{m}}c_{n,k}^{(m)}(f\circ S_{k}\circ T_{n}^{-1})(b)\rightarrow 1\ \mathrm{as}\ m\rightarrow\infty.\]
Hence Lemma \ref{lm2} can be applied. For this, recall that $b$ was assumed to be a smooth point relative to $\Y$, and
note that $(T_{n}\circ S_{k}^{-1}(u),f \circ S_{k}\circ T_{n}^{-1})\in \S_{\Y}\times \S_{\Y^{\ast}}$ are Flinn pairs 
for all $k,n\in \N$, since $(u,f)$ is a Flinn pair by assumption. Thus, by Lemma \ref{lm2} the operator $b^{\ast}\otimes b$ 
can by approximated in the weak operator topology by the finite rank operators
\[P_{m}=\sum_{K_{l_{m}}\times L_{m}} c_{n,k}^{(m)}(f \circ S_{k}\circ T_{n}^{-1})\otimes (T_{n}\circ S_{k}^{-1}(u)).\]
Moreover, Lemma \ref{lm2} implies that b is a Flinn element in $\Y$. A similar argument as in Theorem \ref{bpt} finishes the proof.
\end{proof}

\subsection*{Acknowledgements}
This work is part of the author's Ph.D. research, which is supervised by H.-O. Tylli. The research has been supported 
during 2003-2005 by the Academy of Finland and during 2006 by the Finnish Cultural Foundation.


\begin{thebibliography}{99}
%
%
\bibitem{Am}
D. Amir, \emph{Characterizations of inner product spaces}, Operator Theory: Advances and Applications, vol 20 (1986).       
%
\bibitem{Ba}
S. Banach, \emph{Th\'{e}orie des Op\'{e}rations Lin\'{e}aires}, Warsaw (1932).
%
\bibitem{BR1}
J. Becerra Guerrero \and A. Rodriguez Palacios, The geometry of convex-transitive Banach spaces, 
{\em Bull. London Math. Soc. }31 (1999) 323-331.
%
\bibitem{BR2}
J. Becerra Guerrero \and A. Rodriguez-Palacios, Transitivity of the Norm on Banach Spaces, 
{\em Extracta Math.}17, (2002) 1-58.
%
\bibitem{BR3}
J. Becerra Guerrero \and A. Rodriguez Palacios, Convex-transitive spaces, Big Points, and the Duality Mapping, 
{\em Quart. J. Math. }53, (2002) 257-264.
%
\bibitem{Ca1}
F. Cabello, \emph{10 Variaciones Sobre un Tema de Mazur}, 
Doctoral Thesis, Universidad de Extremadura, (1996).
%
\bibitem{Ca4}
F. Cabello, Transitivity of $M$-spaces and Wood's conjecture,
{\em Math. Proc. Cambridge Phil. Soc. }124 (1998), 513-520.  
%
\bibitem{En}
H.B. Enderton, \emph{Elements of Set Theory}, Academic Press, (1977).
%
\bibitem{HHZ}
P. Habala, P, Hajek, \and V. Zizler, \emph{Introduction to Banach Spaces}, I., matfyzpress, (1996). 
%
\bibitem{KR}
N. Kalton \and B. Randrianantoanina, Surjective isometries of rearrangement invariant subspaces, 
{\em Quart. J. Math. Oxford }45 (1994) 301-327.
%
\bibitem{LTI}
J. Lindenstrauss \and L. Tzafriri, \emph{Classical Banach spaces}. I. Sequence Spaces, Lecture Notes in Mathematics, Vol. 338. 
Springer-Verlag, Berlin-New York, 1973.
%
\bibitem{LTII}
J. Lindenstrauss \and L. Tzafriri, \emph{Classical Banach spaces} II. Function Spaces, Lecture Notes in Mathematics, 
Vol. 338. Springer-Verlag, Berlin-New York, 1973.
%
\bibitem{Ra1}
B. Randrianantoanina, Contractive projections in nonatomic function spaces, 
{\em Proc. Amer. Math. Soc.} 123, (1995).
%
\bibitem{Ra2}
B. Randrianantoanina, A Note on the Banach-Mazur Problem, 
{\em Glasgow Jour. Math. }44, (2002) 159-165.
%
\bibitem{Rol}
S. Rolewicz, \emph{Metric Linear Spaces}, Reidel, Dordrecht, (1985).
%
\bibitem{SZ}
A. Skorik \and M. Zaidenberg, On isometric reflections in Banach spaces, 
{\em Math. Physics, Analysis, Geometry }4 (1997) 212-247. 
%
\bibitem{asypri}
J. Talponen, Asymptotically transitive Banach spaces,
to appear.
%
\bibitem{IR}
J. Talponen, Banach spaces containing isometric reflection subspaces,
in preparation.
%
\bibitem{ZZ}
Z, Zhang \and C. Zhang, On very rotund Banach space, 
{\em Applied Mathematics and Mechanics}, Engl. ed. Vol. 21, (2000).
\end{thebibliography}
\end{document}